\documentclass[12 pt,reqno]{amsart}
\usepackage{amssymb}
\usepackage{amscd}

\newtheorem{lemma}{Lemma}
\newtheorem{prop}[lemma]{Proposition}
\newtheorem{theo}[lemma]{Theorem}
\newtheorem{Kor}[lemma]{Corollary}

\theoremstyle{definition}
\newtheorem{defin}[lemma]{Definition}

\theoremstyle{remark}
\newtheorem{remark}{Remark}
\newtheorem{example}{Example}

\newcommand{\bl}{\bigl<}

\newcommand{\br}{\bigr>}

\newcommand{\G}{\Gamma_g}

\begin{document}
\normalsize
\title[Geometric Representation Theory and $G$-signature]{Geometric Representation
Theory and $G$-signature}
\author{Ralph Grieder}
\address{Department of Mathematics \\ Northwestern University \\ Evanston, 
IL 60208}
\email{ralph@math.nwu.edu}
\thanks{The author is a post-doctoral fellow at Northwestern
University and is
supported by the Swiss National Science Foundation.}

%\date{\today}

\begin{abstract}
Let $G$ be a finite group. To every smooth $G$-action on a
compact, connected and oriented surface we can associate its data of
singular orbits. The set of such data becomes an Abelian group
$\mathbb{B}_G$ under the $G$-equivariant connected sum. We will show
that the map which sends $G$ to $\mathbb{B}_G$ is functorial and
carries many features of the representation theory of finite
groups. We will prove that $\mathbb{B}_G$ consists only of copies of
$\mathbb{Z}$ and $\mathbb{Z}/2\mathbb{Z}$. Furthermore we will show
that there is a surjection from the $G$-equivariant cobordism group of
surface diffeomorphisms to $\mathbb{B}_G$. 

We will define a $G$-signature $\theta$ which is related to the $G$-signature
of Atiyah and Singer and prove that $\theta$ is injective
on the copies of $\mathbb{Z}$ in $\mathbb{B}_G$.
\end{abstract}

\maketitle
\noindent
{\footnotesize {\bf Mathematics Subject Classification (1991).} Primary
57S25, 20C15, 57R85;
Secondary 20F38, 58G10, 11R29.}\\[.3cm]
\noindent
{\footnotesize {\bf Keywords.} Riemann surface, finite group, group action,
equivariant cobordism, G-signature, mapping class group, class number.}

\section{Introduction}\label{section:Intro}
Let $G$ be a finite group and $S_g$ a connected, oriented, compact, Riemann
surface of genus $g$. In these notes we will study diffeomorphic
$G$-actions on surfaces $S_g$ up to isotopy, i.e., we will study all
possible embeddings $\phi$ of the group $G$ into the mapping class
group $\G$. The mapping class group is defined as follows. We write 
$\text{\it Diffeo}_+(S_g)$ for the group of orientation preserving
diffeomorphisms of $S_g$ with the $C^{\infty}$-topology, and
$\text{\it Diffeo}_0(S_g)$ for the connected component of the identity. Then we
have $\G=\text{\it Diffeo}_+(S_g)/\text{\it Diffeo}_0(S_g)$. 

The mapping class group acts on the first homology group of
$S_g$. This action preserves the intersection pairing and thus gives
rise to a symplectic representation $\eta:\G\rightarrow
Sp_{2g}(\mathbb{R})$. The unitary group $U(g)$ is a maximal compact
subgroup of $Sp_{2g}(\mathbb{R})$. Thus for any embedding
$\phi:G\rightarrow\G$ the map $\eta\circ\phi$ factors up to conjugation through
$U(g)\subset Sp_{2g}(\mathbb{R})$ and this unitary representation is denoted
by $\varphi(\phi)$. Hence $\varphi$ assigns to every embedding $\phi$ a
complex representation.

An important and interesting question is to determine the image of the
map $\varphi$. The interest stems, among others, from the fact that
one can use $\varphi$ to deduce results about the cohomology of the
mapping class group. 

Using these ideas in the case where $G$ is the cyclic group of prime
order $p$, it is shown in \cite{Gr1}
that the images of the symplectic classes $d_i\in H^{2i}
(BSp(\mathbb{R});\mathbb{Z})$ in $H^{2i}(\Gamma;\mathbb{Z})$ have
infinite order. Here we consider the stable situation, i.e., $\Gamma$
is the stable mapping class group. 

Another result in \cite{Gr2} states that one can embed polynomial
algebras in the cohomology of the stable mapping class group. For $p$
a regular prime, we have
$H^*(BSp(\mathbb{R});\mathbb{F}_p)\cong\mathbb{F}_p[d_1,d_2,\ldots]$, 
$deg\,d_i=2i$ and the map
$\eta^*:H^*(BSp(\mathbb{R});\mathbb{F}_p)\rightarrow H^*(\Gamma;\mathbb{F}_p)$ is injective on the polynomial algebra
$\mathbb{F}_p[(d_i)_{i\in J}]$, $J=\{i\in\mathbb{N}\,|\,i\equiv 1 \mod{2}\;
\text{or}\; i\equiv 0\mod{p-1}\}$. 

In this paper we want to study the image of the map $\varphi$. We will do
this by turning $\varphi$ into a group homomorphism. First we have to
define a group structure on the set of diffeomorphic $G$-actions. 

We will associate to every $G$-action its singular orbit
data and add two such data by the $G$-equivariant connected sum (see
section \ref{geometry} for definitions). The resulting group of
singular orbit data will be
denoted by $\mathbb{B}_G$. In corollary \ref{ZZ/2Z} we will
prove that the group $\mathbb{B}_G$ consists only of copies of $\mathbb{Z}$ and
$\mathbb{Z}/2\mathbb{Z}$. We will see that the correspondence
$G\mapsto\mathbb{B}_G$ is functorial and that it carries many features of the
representation theory of finite groups. It has induction and
restriction maps and also a double coset formula, therefore we could
describe this functor as a geometric representation theory. 

In section \ref{representation} we will define the map
$\theta:\mathbb{B}_G\rightarrow R_{\mathbb{C}}G/D_G$, which is induced by
$\varphi$. $R_{\mathbb{C}}G/D_G$ denotes a quotient of the complex
representation ring (we
will only consider the additive structure of the representation
ring). We will prove results about the map $\theta$ and its image,
recovering those obtained by Ewing \cite{Ew1} and Edmonds and Ewing
\cite{EdEw} in the special case when $G$ is a finite cyclic group. The
main new result of this section is theorem \ref{Zinj}, which states that
$\theta$ is injective on the copies of $\mathbb{Z}$ in
$\mathbb{B}_G$ for any finite group $G$. 

In section \ref{case} we will apply all the properties of the functor
$\mathbb{B}$ and the map $\theta$ which have been introduced in the
previous sections to the case
$\mathbb{Z}/p\mathbb{Z}\times\mathbb{Z}/p\mathbb{Z}$, $p$ an odd prime. We
will describe the image of $\theta$ by computing its index $\Delta$ in
a subgroup $\mathbb{A}_G$ of the quotient $R_{\mathbb{C}}G/D_G$. This
index will have the form $\Delta=(h^-_p)^{p+1}\cdot p^i$, where
$h^-_p$ is the first factor of the class number of $\mathbb{Q}(e^{2\pi
i/p})$ and $1-p+k\leq i\leq k+(p-1)^2/2$ for some $k\in\mathbb{N}$.

In section \ref{cobordism} we will relate all the results we have
obtained in the previous sections to $G$-equivariant cobordism. We
will construct a surjection $\chi$ from the oriented
$G$-equivariant cobordism group of surface diffeomorphisms
$\Omega_G$ to $\mathbb{B}_G$. With this surjection $\theta$ becomes a
$G$-signature, in the sense that $\theta\circ\chi$ is a homomorphism
from $\Omega_G$ to a quotient of the representation ring of $G$.

It is section \ref{cobordism} which relates this paper to the papers
of Ewing \cite{Ew1} and Edmonds and Ewing \cite{EdEw}. To prove his
results, Ewing \cite{Ew1} essentially looks at the ring
homomorphism from the 
$\mathbb{Z}/p\mathbb{Z}$-equivariant cobordism ring, $p$ a prime, to
the complex representation ring given by the G-signature of Atiyah and Singer
\cite{AtSi} (AS-$G$-signature). 
In the second paper the authors use the AS-$G$-signature as a group
homomorphism from the $\mathbb{Z}/n\mathbb{Z}$-equivariant cobordism
group of surface diffeomorphisms to the complex representation ring.
The value of the AS-$G$-signature for any finite group is computed by the
Atiyah-Bott fixed point theorem \cite{AtBo} and it turns out that it is just
$\varphi(\phi)-\overline{\varphi(\phi)}$; $\overline{\varphi(\phi)}$
denotes the complex conjugate representation. 

We see that the $G$-signature given by $\theta$ and the
AS-$G$-signature are closely related. In the special case of a finite
cyclic group $G$, the results obtained 
in this paper are therefore similar to the ones obtained by Ewing
and Edmonds in \cite{Ew1} and \cite{EdEw}. 

The AS-$G$-signature is very useful as long as the representations remain
complex and non-real (i.e., almost all Abelian
groups, see corollaries \ref{Gabelian} and \ref{repconj}) but for real
representations the AS-$G$-signature is zero. 
The approach described in this paper gives the possibility to
construct a $G$-signature $\theta':\mathbb{B}_G\rightarrow
R_{\mathbb{C}}G/D'_G$ which is injective, even when $\varphi(\phi)$
is a real representation (see remark \ref{S_3}). 

I would like to thank Eric Friedlander for making my stay at
Northwestern University possible. I am also indebted to Stewart Priddy
for his interest and to Shmuel Weinberger for pointing out the
relationship with $G$-equivariant cobordism.

I am especially grateful to Henry Glover for introducing me to this
subject and for his continuing interest and support.

\section{The Group of Singular Orbit Data: Geometry}\label{geometry}
Let $\phi$ be an embedding of any finite group $G$ into some mapping
class group $\G$. By Kerckhoff \cite{Ke} any such embedding can be
lifted to a homomorphism $G\rightarrow\text{\it Diffeo}_+(S_g)$, i.e.,
to an action of $G$ on the surface $S_g$ such that the elements of $G$
act by orientation preserving diffeomorphisms. Let $Gx$ denote the
orbit of $x\in S_g$ under the action of $G$. The orbit is called
singular if $|Gx|<|G|$, else
regular. If the orbit is singular, then there are elements of $G$
which stabilize the point $x\in S_g$ and $G_x$ denotes the
stabilizer of $G$ at $x$. It is proven by Accola \cite[Lemma 4.10]{Ac}
that the stabilizers are cyclic subgroups of $G$. Let $y$ be another element of
the orbit $Gx$. Thus there is an element $a\in G$ such that $ax=y$ and
the stabilizer of $y$ is conjugate to $G_x$, i.e., $aG_xa^{-1}=G_y$. 
As the elements of $G$ operate by orientation preserving
diffeomorphisms there are only finitely many singular orbits. Let
$x_i\in S_g$, $i=1,...,q$ be representatives of these orbits and $\nu_i$
the orders of the stabilizer groups $G_{x_i}$. Every group
$G_{x_i}$ has a generator $\gamma_i\in G$ such that $\gamma_i$ acts by
rotation through $2\pi/\nu_i$ on the tangent space at $x_i$. Similarly
$a\gamma_ia^{-1}$, $a\in G$, generates $aG_{x_i}a^{-1}$ and acts {\it
also} by rotation through $2\pi/\nu_i$ on the tangent space at $ax_i$.
Thus in order to collect information about the singular orbits, it is 
enough to consider the conjugacy classes of the elements $\gamma_i$.
Let $\hat{\gamma_i}$ denote the
conjugacy class of $\gamma_i$ in $G$. The extended singular orbit data
(or extended data) of the embedding $\phi$ is then the collection 
$$\sigma(\phi)=\bigl\{g;\hat{\gamma_1},\ldots,\hat{\gamma_q}\bigr\}_G$$
where $g$ is the genus of the surface $S_g$ and $q$ is the number of
singular orbits of the $G$-action; the conjugacy classes
$\hat{\gamma_i}$ are unique up to order. This data depends a priori on
the chosen lifting of $\phi$, but we will show in the next lemma that
the extended data is well defined for an embedding of
$G$. 
\begin{lemma}
The extended data doesn't depend on the chosen lifting
of the embedding and thus is well defined for any finite subgroup of
the mapping class group.
\end{lemma}
\begin{proof}
Symonds \cite{Sy} proved that the extended data of a diffeomorphism of finite
order depends only upon its isotopy class, and thus is well defined
for an element of finite order of the mapping class group $\G$. But if
the extended data is given for the restriction to any element of a finite
group, then it is also given for the whole group.
\end{proof}
The smaller collection $\bigl\{\hat{\gamma_1},\ldots,\hat{\gamma_q}\bigr\}_G$ 
will be called singular orbit data of the finite group $G$. Notice
that there are infinitely many embeddings of $G$ which give rise to
the same singular orbit data
$\bigl\{\hat{\gamma_1},\ldots,\hat{\gamma_q}\bigr\}_G$. In the sequel
we will omit the subscript $_G$ if it is clear with respect
to which group the conjugacy classes are taken.

We will now define
an addition, the $G$-equivariant connected sum, on the set of
singular orbit data $W_G$ of a group $G$.  

Let $G$ act on a surface $S_g$ with singular orbit data
$\bigl\{\hat{\gamma_1},\ldots,\hat{\gamma_q}\bigr\}$ and on a surface
$S_h$ with singular orbit data
$\bigl\{\hat{\beta_1},\ldots,\hat{\beta_n}\bigr\}$. Find discs $D_1$
in $S_g$ and $D_2$ in $S_h$ such that $\bigl\{aD_j\bigr\}_{a\in G}$
are mutually disjoint for $j=1,2$. Then excise all discs
$\bigl\{aD_j\bigr\}_{a\in G}$, $j=1,2$, from $S_g$ and $S_h$ and take a
connected sum by matching $\partial(aD_1)$ to $\partial(aD_2)$ for all
$a\in G$. The resulting surface $S_{g+h+|G|-1}$ has $|G|$ tubes
joining $S_g$ and $S_h$. The actions of $G$ on $S_g$ and $S_h$ can be
extended to an action on $S_{g+h+|G|-1}$ by permuting the tubes. The
new action has a singular orbit data $\bigl\{\hat{\gamma_1},\ldots,
\hat{\gamma_q},\hat{\beta_1},\ldots,\hat{\beta_n}\bigr\}$.
This construction on surfaces defines a monoid structure on the set
$W_G$ of singular orbit data with addition 
$$\bigl\{\hat{\gamma_1},\ldots,\hat{\gamma_q}\bigr\}\oplus
\bigl\{\hat{\beta_1},\ldots,\hat{\beta_n}\bigr\}:=\bigl\{\hat{\gamma_1},
\ldots,\hat{\gamma_q},\hat{\beta_1},\ldots,\hat{\beta_n}
\bigr\}$$
and zero element $\{\varnothing\}$, the free action. Next we will introduce a
relation to turn this monoid into a group. 

Suppose we have an action of $G$ on a surface $S_g$ with singular orbit data
$\bigl\{\hat{\gamma_1},\hat{\gamma_1}^{-1},\hat{\gamma_2},\ldots,
\hat{\gamma_q}\bigr\}$. Let $\nu=|\bigl<\gamma_1\bigr>|$, then the
conjugacy class $\hat{\gamma_1}$ gives rise to a singular orbit with
representative $x$ such that $a\gamma_1a^{-1}$ acts by rotation through
$2\pi/\nu$ on the tangent space at $ax$. On the other hand
$\hat{\gamma_1}^{-1}$ gives rise to another singular orbit with
representative $z$ such that $a\gamma_1a^{-1}$ acts by rotation through
$-2\pi/\nu$ on the tangent space at $az$. Let $T$ be a set of
representatives for the $\bigl<\gamma_1\bigr>$ left cosets of $G$. Find discs
$D_1$ and $D_2$ around $x$ and $z$ respectively such that $D_j$ is fixed by
$\bigl<\gamma_1\bigr>$, $j=1,2$, and $\cup_{j=1,2}\cup_{t\in
T}\bigl\{tD_j\bigr\}$ are
mutually disjoint. Then excise all discs $\bigl\{tD_j\bigr\}_{t\in
T}$, $j=1,2$, from $S_g$ and connect the boundaries $\partial(tD_1)$
with $\partial(tD_2)$ by means of tubes $S^1\times[0,1]$ for every
$t\in T$. The resulting surface $S_{g+w}$ has $w=|G|/\nu$ new
handles. The action of $G$ on $S_g$ can be extended to $S_{g+w}$ by
permuting and rotating the new handles. This extended action yields
the singular orbit data
$\bigl\{\hat{\gamma_2},\ldots,\hat{\gamma_q}\bigr\}$. Pairs
of singular orbits which have opposite rotation on the tangent spaces
will be called cancelling pairs. The above process of eliminating
cancelling pairs will be called reduction and if there are no such
cancelling pairs left the singular orbit data is said to be in reduced form. 

Now we can define the relation.
$$\bigl\{\hat{\gamma_1},\ldots,\hat{\gamma_q}\bigr\}\sim
\bigl\{\hat{\beta_1},\ldots,\hat{\beta_n}\bigr\}:\Leftrightarrow
\left\{\begin{array}{l}
\text{The two singular orbit data}\\
\text{have the same reduced form}
\end{array}
\right.$$
This relation defines an equivalence relation on the set $W_G$ of
singular orbit data. 
$$\mathbb{W}_G:=W_G/\sim$$
The set $\mathbb{W}_G$ is not only a commutative monoid as $W_G$ but
contains also inverse elements and thus is a commutative group. The
inverse element of
$\bigl\{\hat{\gamma_1},\ldots,\hat{\gamma_q}\bigr\}$ is
$\bigl\{\hat{\gamma_1}^{-1},\ldots,\hat{\gamma_q}^{-1}\bigr\}$ and the
zero elements are cancelling pairs
$\bigl\{\hat{\gamma},\hat{\gamma}^{-1}\bigr\}$, 
$\gamma\in G$, and $\{\varnothing\}$ the free action.

\section{The Group of Singular Orbit Data: Algebra}
In the last section we defined a group structure on the singular orbit
data. In this section we will define the same group but in an
algebraic context. We will also find necessary and sufficient
conditions for a $q$-tuple of conjugacy classes to define a singular
orbit.

Every embedding $G\hookrightarrow\G$ with extended data
$\bigl\{g;\hat{\gamma_1},\ldots,\hat{\gamma_q}\bigr\}_G$
gives rise to a surjective homomorphism
\begin{equation}\label{p}
\bl a_1,b_1,\ldots,a_h,b_h,\xi_1,\ldots,\xi_q\,|\,
\prod_{i=1}^h[a_i,b_i]\cdot\xi_1\cdots\xi_q\br\stackrel{p}{\twoheadrightarrow}G
\end{equation}
such that the Riemann-Hurwitz equation 
\begin{equation}\label{RiHu}
2g-2=|G|(2h-2)+|G|\Bigl(\sum_{i=1}^q1-1/\nu_i\Bigr)
\end{equation}
is satisfied and $p(\xi_i)=\gamma_i$, $i=1,\ldots,q$.
Here $h$ denotes the genus of the quotient surface $S_g/G$, $q$ the
number of singular orbits and $\nu_i$ the order of the cyclic subgroup
generated by $p(\xi_i)$.

By the discussion in Accola's book \cite[Section 4.10 Branched
Coverings]{Ac} the converse is also true, for any surjective
homomorphism $p$ as in (\ref{p}), there is an embedding
$G\hookrightarrow\G$ with extended data
$\bigl\{g;\hat{\gamma_1},\ldots,\hat{\gamma_q}\bigr\}_G$ where $g$
satisfies the Riemann-Hurwitz equation (\ref{RiHu}). 

Furthermore any surjection $p$ as in (\ref{p}) gives rise to a $q$-tuple
$\bigl(p(\xi_1),..$ $..,p(\xi_q)\bigr)$ such that $p(\xi_1)\cdots p(\xi_q)\in
[G,G]$, the commutator subgroup of $G$. 

Conversely for any $q$-tuple
$\big(\gamma_1,\ldots,\gamma_q\big)\in G\times\cdots\times G$, with
$\gamma_1\cdots \gamma_q\in [G,G]$ one can construct a surjection $p$
as follows.
Map $\xi_i$ to $\gamma_i$ for $i=1,\ldots,q$ and choose
appropriate images for the $a_i$'s and $b_i$'s in $G$ such that $p$
becomes surjective and the relation
$\prod_{i=1}^h[a_i,b_i]\cdot\xi_1\cdots\xi_q$ is satisfied. This is
always possible as one can choose $h$ big enough. Observe
that the construction of $p$ doesn't depend on the ordering of the
$\gamma_i$ and neither on the representative of the conjugacy class of
$\gamma_i$. 

The above discussion shows that there is a correspondence between a
singular orbit data
$\bigl\{\hat{\gamma_1},\ldots,\hat{\gamma_q}\bigr\}_G$ and an
unordered $q$-tuple of conjugacy classes
$\bigl(\hat{\gamma_1},\ldots,\hat{\gamma_q}\bigr)$ with
$\gamma_1\cdots \gamma_q\in [G,G]$. Later we are going to make
the correspondence more precise. But first consider the following map.
\begin{eqnarray}
\Psi_q:G\times\cdots\times G&\longrightarrow&\frac{G}{[G,G]} \\
\big(\gamma_1,\ldots,\gamma_q\big)&\longmapsto&
\bar{\gamma_1}\cdots\bar{\gamma_q}\notag
\end{eqnarray}
$\Psi_q$ is a homomorphism and every element of $ker\,\Psi_q$ gives
rise to a map $p$ by the above considerations. It is clear by the
definition of $\Psi_q$ that $ker\,\Psi_q$ is invariant under
permutation and conjugation of the components. 
\begin{defin}
$$\Lambda_G:=\bigcup_{q=0}^{\infty}ker\,\Psi_q/\sim$$
Here $\sim$ means up to permutation and conjugation of the
components and $ker\,\Psi_0$ is defined to be the empty tuple. 
An element of $\Lambda_G$ will be denoted by 
$\bigl[\hat{\gamma_1},\ldots,\hat{\gamma_q}\bigr]_G$, i.e., an unordered
$q$-tuple of conjugacy classes of $G$, with
$\gamma_1\cdots\gamma_q\in[G,G]$. We will consider $ker\,\Psi_q/\sim$ as a
subset of $ker\,\Psi_{q+1}/\sim$ by the inclusion
$\bigl[\hat{\gamma_1},\ldots,\hat{\gamma_q}\bigr]_G\mapsto\bigl[\hat{\gamma_1},
\ldots,\hat{\gamma_q},1\bigr]_G$, in particular all the q-tuples of ones
$\bigl[1,\ldots,1\bigr]_G$ will be the same as the empty tuple
$[\varnothing]_G$.
\end{defin}
One can now define a commutative monoid structure on $\Lambda_G$ by 
$$\bigl[\hat{\gamma_1},\ldots,\hat{\gamma_q}\bigr]_G\oplus
\bigl[\hat{\beta_1},\ldots,\hat{\beta_r}\bigr]_G:=\bigl[\hat{\gamma_1},
\ldots,\hat{\gamma_q},\hat{\beta_1},\ldots,\hat{\beta_r}\bigr]_G.$$
The associativity and commutativity of this addition is immediate, to
obtain a group structure one has to introduce inverse elements. This
is done by the next definition.
\begin{defin}
$$\mathbb{B}_G:=\Lambda_G/\bl[\hat{\gamma},\hat{\gamma}^{-1}]_G\,,\,\gamma\in
G\br$$
\end{defin}
$\mathbb{B}_G$ is now a commutative group. The zero elements are 
$[\hat{\gamma},\hat{\gamma}^{-1}]_G=[\varnothing]_G=0$ and the inverse
elements $\ominus
\bigl[\hat{\gamma_1},\ldots,\hat{\gamma_q}\bigr]_G=\bigl[\hat{\gamma_1}^{-1},\ldots,\hat{\gamma_q}^{-1}\bigr]_G$.
Now we can be more precise about the correspondence between the
singular orbit data and the elements of $\mathbb{B}_G$. 
\begin{prop}
There is an isomorphism of groups
\begin{eqnarray*}
\mathbb{B}_G&\stackrel{\sim}{\rightarrow}&\mathbb{W}_G\\
\bigl[\hat{\gamma_1},\ldots,\hat{\gamma_q}\bigr]_G&\mapsto&
\bigl\{\hat{\gamma_1},\ldots,\hat{\gamma_q}\bigr\}_G
\end{eqnarray*}
given by the correspondence:\\
An element $\bigl[\hat{\gamma_1},\ldots,\hat{\gamma_q}\bigr]_G$ of
$\mathbb{B}_G$ gives rise to a map $p$ as in (\ref{p}) which in turn
defines an embedding $G\hookrightarrow\G$ with singular orbit data
$\bigl\{\hat{\gamma_1},\ldots,\hat{\gamma_q}\bigr\}_G$.
\end{prop}
\begin{proof}
The map is well defined as the different choices of $p$ give rise to
the same singular orbit data. 

By the construction it follows
immediately that the map is a homomorphism and injective. 

We can invert the correspondence by assigning to every singular orbit data
$\bigl\{\hat{\gamma_1},\ldots,\hat{\gamma_q}\bigr\}_G$ a map $p$ as in
(\ref{p}) and thus the element
$\bigl[\hat{\gamma_1},\ldots,\hat{\gamma_q}\bigr]_G$ in
$\mathbb{B}_G$. This shows that the map is surjective.
It follows also that the correspondence is one-to-one.
\end{proof}
From now on we will refer to $\mathbb{B}_G$ as the group of singular
orbit data and to elements
$\bigl[\hat{\gamma_1},\ldots,\hat{\gamma_q}\bigr]_G$ as singular orbit
data. If it is clear which group we mean, we will omit the
subscript $_G$ and simply write
$\bigl[\hat{\gamma_1},\ldots,\hat{\gamma_q}\bigr]$ for an element of
$\mathbb{B}_G$. 

For every group homomorphism $f:H\rightarrow G$ there is an induced
homomorphism of their associated groups of singular orbit data.
\begin{eqnarray*}
\mathbb{B}_f:\mathbb{B}_H & \rightarrow & \mathbb{B}_G \\
\bigl[\hat{\gamma_1},\ldots,\hat{\gamma_q}\bigr]_H & \mapsto &
\bigl[\widehat{f(\gamma_1)},\ldots,\widehat{f(\gamma_q)}\bigr]_G
\end{eqnarray*}
This map is well defined as $f(\gamma_1)\cdots
f(\gamma_q)=f(\gamma_1\cdots\gamma_q)\in f([H,H])\subset [G,G]$. It is
also true that $\mathbb{B}_{f\circ g}=\mathbb{B}_f\circ\mathbb{B}_g$ thus
we have a covariant functor $\mathbb{B}$ from the category of finite
groups into the category of Abelian groups. 
$$\mathbb{B}:G\mapsto\mathbb{B}_G$$
\begin{remark}
If the map $f$ is injective and the groups $G$ and $H$ are Abelian,
then the map $\mathbb{B}_f$ is also injective.
\end{remark}
\begin{remark}
If the map $f$ is an isomorphism then $\mathbb{B}_f$ is also an
isomorphism. In particular conjugation by an element of the group $G$ will
induce the identity map on $\mathbb{B}_G$, thus we have a map
$Out(G)\rightarrow Aut(\mathbb{B}_G)$.
\end{remark}
\begin{remark}
For every integer $n$ we have the following action
$$n*\bigl[\hat{\gamma_1},\ldots,\hat{\gamma_q}\bigr]:=\bigl[\hat{\gamma_1^n},\ldots,\hat{\gamma_q^n}\bigr].$$
This defines a map $\mathbb{Z}\stackrel{\lambda}{\rightarrow}End(\mathbb{B}_G)$
which preserves the multiplicative structure of $\mathbb{Z}$,
$\lambda(nm)=\lambda(n)\circ\lambda(m)$, in particular $\lambda(0)$ is
the zero map and $\lambda(1)$ is the identity map. If $m$ is the order
of $G$ then the map $\lambda$ factors through
$\mathbb{Z}/m\mathbb{Z}$. The group of singular orbit data with its
$\mathbb{Z}$-action is an invariant of the group $G$. 
\end{remark}
In the case where $K$ is a subgroup of $G$ then there is another map
$\mathbb{B}res^G_K:\mathbb{B}_G\rightarrow \mathbb{B}_K$ called the restriction
map. It is defined by restricting the singular orbits of $G$ to the
subgroup $K$. If the group $G$ is Abelian we can give an explicit
formula for this map, in the general case this formula becomes too
complicated. We will omit the $\hat{\;}$ notation for Abelian groups.
$$\mathbb{B}res^G_K\bigl(\bigl[\gamma_1,\ldots,\gamma_q\bigr]\bigr)=\bigl[\gamma_{11}^
{n_1},\ldots,
\gamma_{1i_1}^{n_1},\ldots,\gamma_{q1}^{n_q},\ldots,\gamma_{qi_q}^{n_q}\bigr]$$
with $i_r=\frac{|G|\cdot|K\cap<\gamma_r>|}{|<\gamma_r>|\cdot |K|}$,
$n_r=\text{min}\{m|\gamma_r^m\in K\}$ and $\gamma_{rs}=\gamma_r$ for
$s=1,\ldots,i_r$. 

Let $i:H\rightarrow G$ be another subgroup, then there exists a
double coset formula for singular orbit data.
\begin{prop}\label{doublecoset}
Let $H$ and $K$ be subgroups of $G$ and $S$ a set of representatives
for the $(H,K)$ double cosets of $G$. For $s\in S$, let
$H_s=sHs^{-1}\cap K$, which is a subgroup of $K$. The inclusions
will be denoted by $i:H\rightarrow G$ and $j:H_s\rightarrow K$. Then
we have the following equation
$$\mathbb{B}res^G_K\circ\mathbb{B}_i= 
\bigoplus_{s\in S}\mathbb{B}_j\circ
\mathbb{B}res^{sHs^{-1}}_{H_s}\circ\mathbb{B}_{f_s}$$
where $f_s$ is just conjugation by $s$, i.e., $f_s(g)=sgs^{-1}$.
\end{prop}
\begin{proof}
Let
$\bigl[\hat{\gamma},\ldots\bigr]_H$ be a singular orbit data in
$\mathbb{B}_H$, then we have $\mathbb{B}_i(\bigl[\hat{\gamma},\ldots\bigr]_H)
=\bigl[\hat{\gamma},\ldots\bigr]_G$. 
First we consider the restriction on the level of singular orbits. Let
$Gx$ be the singular orbit corresponding to $\hat{\gamma}$. Thus
$\gamma$ operates by rotation through $2\pi/|\bigl<\gamma\bigr>|$ on $x$. We can
now partition the singular orbit $Gx$. 
$$Gx=\bigcup_{s\in S}\bigcup_{g\in KsH}gx=\bigcup_{s\in S}\bigcup_{h\in H}
Kshx=\bigcup_{s\in S}\bigcup_{h\in sHs^{-1}/H_s}Khsx$$
Thus we have 
$$Gx=\bigcup_{s\in S}\bigcup_{h\in sHs^{-1}/H_s}Khsx$$
Note that the unions in the last equation are disjoint.
Now we know how the restricted singular orbits look like. But what are
the stabilizers of $K$ at $hsx$?

$K\cap\bigl<hs\gamma s^{-1}h^{-1}\bigr>=s^{-1}h^{-1}Khs\cap\bigl<\gamma\bigr>=\bigl<\gamma_{_{hs}}\bigr>$
and $\gamma_{_{hs}}$ is unique such that it operates by rotation through
$2\pi/|\bigl<\gamma_{_{hs}}\bigr>|$ on $x$. Thus we conclude:
\begin{eqnarray*}
\mathbb{B}res^G_K\circ\mathbb{B}_i(\bigl[\hat{\gamma},\ldots\bigr]_H)&=&
\mathbb{B}res^G_K(\bigl[\hat{\gamma},\ldots\bigr]_G) \\
&=&\bigl[\widehat{hs\gamma_{_{hs}}
(hs)^{-1}},\ldots\bigr]_K\;,\;h\in sHs^{-1}/H_s\;,\;s\in S
\end{eqnarray*}
On the other hand we have
\begin{eqnarray*}
\mathbb{B}_j\circ \mathbb{B}res^{sHs^{-1}}_{H_s}\circ\mathbb{B}_{f_s}
\bigl[\hat{\gamma},\ldots\bigr]_H &=&\mathbb{B}_j\circ \mathbb{B}res^{sHs^{-1}}_{H_s}
\bigl[\widehat{s\gamma s^{-1}},\ldots\bigr]_{sHs^{-1}}= \\
\mathbb{B}_j\bigl[\widehat{hs\gamma_{_{hs}}(hs)^{-1}},\ldots\bigr]_{H_s}& =&
\bigl[\widehat{hs\gamma_{_{hs}}(hs)^{-1}},\ldots\bigr]_K\;\;\;,\;h\in sHs^{-1}/H_s
\end{eqnarray*}
and both sides of the equation in the lemma are the same.
\end{proof}
In this section we showed that many features of the representation
theory of finite groups carry over to the functor $\mathbb{B}$.
Thus we could say that $\mathbb{B}$ describes a geometric representation
theory. 

In the following examples the cyclic group of order $m$ will be
denoted by $C_m$. 
\begin{example}
$G=C_2$; $\mathbb{B}_G\cong {0}$.
\end{example}
\begin{example}\label{exCmodd}
Let $m$ be odd, $G=C_m=\bigl<x\bigr>$ then $\mathbb{B}_G\cong
\mathbb{Z}^{\frac{m-1}{2}}$ and a basis for $\mathbb{B}_G$ is given by
$\bigl[x,x^i,x^{m-i-1}\bigr]$, $i=1,\ldots,\frac{m-1}{2}$.
\end{example}
\begin{example}\label{exCmeven}
Let $m$ be even, $G=C_m=\bigl<x\bigr>$ then $\mathbb{B}_G\cong
\mathbb{Z}^{\frac{m}{2}-1}$ and a basis for $\mathbb{B}_G$ is given by
$\bigl[x,x^i,x^{m-i-1}\bigr]$, $i=1,\ldots,\frac{m}{2}-1$.
\end{example}
\begin{example}\label{CpCp}
Let $p$ be an odd prime, $G=C_p\times C_p=\bigl<x\bigr>\times\bigl<y\bigr>$
then $\mathbb{B}_G\cong\mathbb{Z}^{\frac{p^2-1}{2}}$ and a basis for
$\mathbb{B}_G$ is given by $[x^j,y^i,x^{p-j}y^{p-i}]$,
$j=1,\ldots,$\linebreak
 $(p-1)/2$, $i=1,\ldots,p-1$; $[x,x^k,x^{p-k-1}]$,
$k=1,\ldots,(p-1)/2$; $[y,y^l,y^{p-l-1}]$, $l=1,\ldots,(p-1)/2$.\\[2mm]
For $p$ the even prime, we have $\mathbb{B}_G\cong C_2$ and the
generator is $\bigl[x,y,xy\bigr]$.
\end{example}
\begin{example}\label{exS_3}
Let $G=S_3=C_3\rtimes C_2$, $C_3=\bigl<a\bigr>$ and $C_2=\bigl<b\bigr>$. Then we have
$[G,G]=C_3$, $\hat{a}=\{a,a^2\}$ and $\hat{b}=\{ab,a^2b,b\}$. The
singular orbit data are then $\bigl[\hat{a}\bigr]$,
$2\cdot\bigl[\hat{a}\bigr]=\bigl[\hat{a},\hat{a}\bigr]=\bigl[\hat{a},\hat{a^2}\bigr]=0$,
$\bigl[\hat{b},\hat{b}\bigr]=0$. Thus the group of singular orbit data
is generated by $\bigl[\hat{a}\bigr]$ and
$\mathbb{B}_{S_3}=\bigl<\bigl[\hat{a}\bigr]\bigr>\cong C_2$. 

As $\mathbb{B}_{C_2}\cong {0}$ the only maps which are of any interest are
the maps
$\mathbb{B}res^{S_3}_{C_3}:\mathbb{B}_{S_3}\rightarrow\mathbb{B}_{C_3}$
and $\mathbb{B}_i:\mathbb{B}_{C_3}\rightarrow\mathbb{B}_{S_3}$,
where $i$ is the inclusion $i:C_3\hookrightarrow S_3$. The maps are given by 
$\mathbb{B}res^{S_3}_{C_3}([\hat{a}]_{S_3})=[a,a^2]_{C_3}=0$ and 
$\mathbb{B}_i([a,a,a]_{C_3})=[\hat{a},\hat{a},\hat{a}]_{S_3}=
[\hat{a}]_{S_3}$. Thus $\mathbb{B}res^{S_3}_{C_3}$ is just the zero
map and $\mathbb{B}_i$ is the surjection of $\mathbb{Z}$ onto $C_2$.
\end{example}
\begin{theo}\label{torsion}
The only torsion $\mathbb{B}_G$ contains, $G$ any finite group, is two-torsion.
\end{theo}
\begin{proof}
Let $\alpha=\bigl[\hat{\gamma_1},\ldots,\hat{\gamma_q}\bigr]$ be an element
of $\mathbb{B}_G$. Suppose $\alpha$ doesn't contain any cancelling
pairs and it has finite order, i.e., 
$$n\cdot\bigl[\hat{\gamma_1},\ldots,\hat{\gamma_q}\bigr]=0.$$
It follows that $n\cdot\alpha$ consists of cancelling pairs and thus
$\hat{\gamma_i}=\hat{\gamma_i}^{-1}$ and $n$ has to be even. But then
we have already $2\cdot\alpha=0$, which proves the theorem.
\end{proof}
\begin{Kor}\label{ZZ/2Z}
$$\mathbb{B}_G\cong\mathbb{Z}\times\cdots\times\mathbb{Z}\times\mathbb{Z}/2\mathbb{Z}\times\cdots\times\mathbb{Z}/2\mathbb{Z}$$
The copies of $\mathbb{Z}/2\mathbb{Z}$ are given by elements
$\bigl[\hat{\gamma_1},\ldots,\hat{\gamma_q}\bigr]$ such that when in
reduced form $\gamma_i$ is conjugate to $\gamma_i^{-1}$, $i=1,\ldots,q$.
\end{Kor}
\begin{proof}
$\mathbb{B}_G$ is an Abelian group and by theorem \ref{torsion} it
contains only two-torsion. The group $\mathbb{B}_G$ is also finitely
generated. Indeed, every element of $\mathbb{B}_G$ can be written as a
linear combination of triples
$\bigl[\hat{\gamma_i},\hat{\gamma_j},\hat{\gamma_k}\bigr]$ (this fact
is proven in \cite{Gr3}) and there are only finitely many triples for
a finite group $G$. The first assertion follows now by the
structure theorem for Abelian groups. 

If $\bigl[\hat{\gamma_1},\ldots,\hat{\gamma_q}\bigr]$ is a reduced element of
order two then
$$\bigl[\hat{\gamma_1},\ldots,\hat{\gamma_q}\bigr]=\ominus\bigl[\hat{\gamma_1}
,\ldots,\hat{\gamma_q}\bigr]=\bigl[\hat{\gamma_1}^{-1},\ldots,
\hat{\gamma_q}^{-1}\bigr]$$
and $\gamma_i$ is conjugate to $\gamma_i^{-1}$, $i=1,\ldots,q$. The
converse however is also true.
\end{proof}
\begin{Kor}\label{BGtorsion}
Let $G$ be an Abelian group which doesn't contain a copy of
$\mathbb{Z}/2\mathbb{Z}\times\mathbb{Z}/2\mathbb{Z}$, then there is no
two torsion in $\mathbb{B}_G$.
\end{Kor}
\begin{proof}
Assume that $\bigl[\gamma_1,\ldots,\gamma_q\bigr]\in\mathbb{B}_G$ be a
reduced element of order two. Thus $\gamma_j=\gamma_j^{-1}$,
$j=1,\ldots,q$ and $q>2$, i.e., there are at least two distinct
involutions in $G$ which commute. This contradicts the assumption on $G$.
\end{proof}
In \cite{Gr3} the author gives a complete computation
of the group $\mathbb{B}_G$ for any finite group $G$ .

\section{Representations}\label{representation}
In the previous sections we defined a group structure on the set of isotopy
classes of diffeomorphic group actions. The resulting group was
denoted by $\mathbb{B}_G$. In this section we will define
a group homomorphism $\theta$ from the group $\mathbb{B}_G$ to a quotient of
the complex representation ring of $G$. First we need a proposition
which shows how to compute the 
representation $\varphi(\phi)$ from a given $G$-action $\phi$. The map
$\theta$ will then be induced by $\varphi$. We start with cyclic
groups and follow 
the presentation in Farkas and Kra's book \cite[pp 274]{FaKr}.

\subsection{Cyclic Groups}\label{cyclic}
Let $G$ be the cyclic group $C_n=\bigl<x\bigr>$ and $\phi$ an embedding of
$G$ with extended data $\sigma(\phi)=\bigl\{g;$
$x^{i_1},\ldots,x^{i_q}\bigr\}_{C_n}$. The representation $\varphi(\phi)$
will be a linear combination $\varphi(\phi)=\sum_{j=0}^{n-1}n_j\rho_i$ where
$n_j\in\mathbb{N}$ and the
$\rho_i$'s, $i=1,\ldots,n-1$ are the one dimensional 
irreducible representations, such that $\rho_1$ is a faithful
representation and $\rho_i=\rho_1^{\otimes i}$; $\rho_0$ is the one dimensional
trivial representation. We will give a formula for the multiplicities
$n_j$ but first we need some notations:\\
$gcd(n,i_s)=k_s$, $i_s=l_sk_s$ and thus $l_s$ is prime to $n$. Let
$r_{sj}$ be the unique integer such that
$$1\leq r_{sj}\leq n/k_s\;\;\;\text{and}\;\;\;r_{sj}\equiv jl_s \mod{n/k_s}\;,\;\;j=1,\ldots,n-1.$$
The formula in Farkas and Kra's book \cite[pp 274]{FaKr} reads now:
\begin{eqnarray}\label{n_j}
n_j & = &h-1+q-\big(\sum_{s=1}^qk_sr_{sj}\big)/n\;,\;\;j=1,\ldots,n-1 \\
n_0 & = & h.\notag
\end{eqnarray}
Here $h$ denotes the genus of the quotient surface, it is well defined
by the Riemann-Hurwitz equation (\ref{RiHu}). Thus we have
proven the following lemma.
\begin{lemma}\label{repcyclic}
For a cyclic group $C_n$ and an embedding $\phi$ the representation
$\varphi(\phi)$ can be computed by its extended data
$\sigma(\phi)=\bigl\{g;$ $x^{i_1},\ldots,x^{i_q}\bigr\}_{C_n}$. The formula for the
multiplicities $n_j$ is given by equation (\ref{n_j}).
\end{lemma}
In the sequel the regular representation of $G$ will be denoted by $\omega_G$
or $\omega_m$ if the group is cyclic of order $m$. 
\begin{example}\label{freesing}
Let $G=C_n$. For a free action $\tilde{\phi}$ we have $\sigma(\tilde{\phi})=\bigl\{g;\varnothing\bigr\}_{C_n}$ 
and the Riemann-Hurwitz equation (\ref{RiHu}) reads $g=n(h-1)+1$. Equation
(\ref{n_j}) tells us that $\varphi(\tilde{\phi} )=\rho_0+ (h-1)\cdot\omega_n$. 
\end{example}
\begin{example}\label{zerosing}
In this example we will consider the extended data
$\bigl\{g;x^i,$ $x^{n-i}\bigr\}_{C_n}$. In this situation we have
$gcd(n,i)=k=gcd(n,n-i)$, $i=lk$, $n-i=(n/k-l)k$ thus $r_{1j}\equiv jl
\mod{n/k}$ and $r_{2j}\equiv j(n/k-l)\mod{n/k}$. This implies
$r_{2j}=n/k-r_{1j}$ if $r_{1j}\neq n/k$ and $r_{2j}=r_{1j}$ if
$r_{1j}=n/k$. For the $n_j$ we obtain
\begin{eqnarray*}
n_j&=&h-1+2-(kr_{1j}/n)-(k(n/k-r_{1j})/n)=h\;,\;j\not\equiv
0\mod{n/k} \\
n_j&=&h-1+2-1-1=h-1 \;\;,\;\;j\equiv 0\mod{n/k} \\
n_0&=&h\\
\end{eqnarray*}
with Riemann-Hurwitz equation (\ref{RiHu}) $g=hn+(1-k)$. We can now
express $\varphi(\phi)$ in terms of induced representations.
$$\varphi(\phi)=\sum_{j=0}^{n-1}n_j\rho_i=h\cdot Ind^{C_n}_{<0>}\rho_0+\rho_0-Ind^{C_n}_{<x^i>}\rho_0=h\cdot\omega_n+\rho_0-Ind^{C_n}_{C_{n/k}}\rho_0$$
\end{example}
It follows readily by examples \ref{freesing} and \ref{zerosing} that
the additive subgroup of $R_{\mathbb{C}}C_n$ 
generated by all the representations induced by the extended data
$\bigl\{g;x^i,x^{n-i}\bigr\}_{C_n}$, $i=0,\ldots,n-1$, is the
same as the subgroup 
$$D_{C_n}:=\bigl<Ind^{C_n}_{H}\rho_0\;\bigm|\;H\leq C_n\bigr>.$$
\begin{lemma}\label{integralcyc}
Every character with values in $\mathbb{Z}$ is a linear combination,
with coefficients in $\mathbb{Z}$, of characters
$Ind^{C_n}_{H}\rho_0$, where $H$ runs through all the subgroups
of $C_n$, i.e.,
$\bigl<Ind^{C_n}_{H}\rho_0\;\bigm|\;H\leq C_n\bigr>=R_{\mathbb{Q}}C_n.$
\end{lemma}
\begin{proof}
For each divisor $d$ of $n$, let 
\begin{eqnarray*}
\chi_d&=&\sum_{l=1,gcd(l,d)=1}^{d-1}\rho_l^{\otimes
n/d}\;\;,\;\;d\neq1\\
\chi_1&=&\rho_0\;\;,\;\;d=1
\end{eqnarray*}
Then the $\chi_d$ form an orthogonal basis of $R_{\mathbb{Q}}C_n$. On
the other hand we have 
$$Ind^{C_n}_{C_{n/d}}\rho_0=\sum_{l=1}^{d}\rho_l^{\otimes
n/d}=\sum_{d'|d}\sum_{l=1,gcd(l,d')=1}^{d'-1}\rho_l^{\otimes
n/d'}=\sum_{d'|d}\chi_{d'}.$$
By induction we can now prove that the $\chi_d$ are linear
combinations of the $Ind^{C_n}_{C_{n/d'}}\rho_0$, $d'|d$, with
integer coefficients. Thus the $Ind^{C_n}_{C_{n/d}}\rho_0$, $d|n$,
form a basis of $R_{\mathbb{Q}}C_n$. (See also exercise 13.1 of
Serre's book \cite{Se}.)
\end{proof}
Let $\sigma(\tilde{\phi})$ and $\sigma(\bar{\phi})$ be two extended
data which have the same singular orbits. Then by
equation (\ref{n_j}) their
associated unitary representations $\varphi(\tilde{\phi})$ and
$\varphi(\bar{\phi})$ differ only by a multiple of the regular
representation $\omega_G$. Thus by example \ref{freesing} there is a
well defined map 
$$\Lambda_{C_n}\rightarrow R_{\mathbb{C}}C_n/\bigl<\rho_0,\omega_n\bigr>$$
which is induced by $\varphi$. This map is even a homomorphism of
Abelian monoids. Indeed, given two extended data
$\sigma(\phi)=\bigl\{g;x^{i_1},\ldots,x^{i_q}\bigr\}_{C_n}$ and
$\sigma(\phi')=\bigl\{g';$ $x^{j_1},\ldots,x^{j_w}\bigr\}_{C_n}$. 
Then again by equation (\ref{n_j}) we can see that
$$\varphi(\phi)+\varphi(\phi')\equiv\varphi(\tau)\;\;\;\text{up to
multiples of $\rho_0$ and $\omega_n$,}$$
where $\tau$ gives rise to an extended data 
$\bigl\{\bar{g};x^{i_1},\ldots,x^{i_q},$ 
$x^{j_1},\ldots,x^{j_w}\bigr\}_{C_n}$.

Example \ref{zerosing} shows that the zero-elements of
$\mathbb{B}_{C_n}$ are mapped into $D_{C_n}$ thus we obtain a well
defined homomorphism of Abelian groups which is induced by $\varphi$:
$$\theta: \mathbb{B}_{C_n}\rightarrow
R_{\mathbb{C}}C_n/D_{C_n}=:\widetilde{R_{\mathbb{C}}C_n}.$$
Let $H$ be a subgroup of $C_n$. Then the homomorphisms
$Res^{C_n}_H:R_{\mathbb{C}}C_n\rightarrow R_{\mathbb{C}}H$ and
$Ind^{C_n}_H:R_{\mathbb{C}}H\rightarrow R_{\mathbb{C}}C_n$ map
$D_{C_n}$ into $D_H$ and $D_H$ into $D_{C_n}$ respectively. Thus they
induce well defined homomorphisms
$\widetilde{Res}^{C_n}_H:\widetilde{R_{\mathbb{C}}C_n}\rightarrow
\widetilde{R_{\mathbb{C}}H}$ and
$\widetilde{Ind}^{C_n}_H:\widetilde{R_{\mathbb{C}}H} \rightarrow
\widetilde{R_{\mathbb{C}}C_n}$.
\begin{prop}\label{Ressquare}
For a cyclic group $G$ and a subgroup $H<G$ the following diagram commutes.
\begin{eqnarray*}
\mathbb{B}_H & \stackrel{\mathbb{B}res^G_H}{\leftarrow} & \mathbb{B}_G \\
\downarrow & & \downarrow \\
\widetilde{R_{\mathbb{C}}H} & \stackrel{\widetilde{Res}^G_H}{\leftarrow} & \widetilde{R_{\mathbb{C}}G}
\end{eqnarray*}
\end{prop}
\begin{proof}
Let $\alpha\in\mathbb{B}_G$ and $\phi_{\alpha}:G\hookrightarrow\G$
denote a corresponding embedding. We can now easily describe a
possible embedding for $\mathbb{B}res^G_H\alpha$, namely
$\phi_{\mathbb{B}res_H^G\alpha}=\phi_{\alpha}|_H:H<G\hookrightarrow\G$. From
this we deduce
$Res^G_H\varphi(\phi_{\alpha})=\varphi(\phi_{\mathbb{B}res^G_H\alpha})$
which proves the proposition.
\end{proof}
\begin{prop}\label{Indsquare}
For a cyclic group $G$ and a subgroup $H<G$ the following diagram commutes.
\begin{eqnarray*}
\mathbb{B}_H & \stackrel{\mathbb{B}_i}{\rightarrow} & \mathbb{B}_G \\
\downarrow &  & \downarrow \\
\widetilde{R_{\mathbb{C}}H} & \stackrel{\widetilde{Ind}^G_H}{\rightarrow} & \widetilde{R_{
\mathbb{C}}G}
\end{eqnarray*}
\end{prop}
\begin{proof}
For the proof we need some explicit computations, so we choose
$G=C_n=\bigl<x\bigr>$ and $H=C_l=\bigl<y\bigr>$ where $l|n$ with the
inclusion $i:H\rightarrow G$, $y\mapsto x^{n/l}$. 
In the sequel we will use the notation $r_{sj}^l$, $n_j^l$ and 
$r_{sj}^n$, $n_j^n$ to denote integers dealing with the groups $C_l$
and $C_n$ respectively. 

Let $\alpha=[y^{i_1},\ldots,y^{i_q}]$ be an element of $\mathbb{B}_H$. Then
with the notation of equation (\ref{n_j}) we obtain:
$gcd(l,i_s)=k_s$, $s=1,\ldots,q$, $i_s=k_su_s$ and thus
$$r_{sj}^l\equiv ju_s\mod{l/k_s}\;\;,\;\;j=1,\ldots,l-1$$
where $1\leq r_{sj}^l\leq l/k_s$. By equation (\ref{n_j}) we obtain
for the multiplicities of the representations:
$$-n_j^l\equiv \big(\sum_{s=1}^qk_sr_{sj}^l\big)/l\;,\;\text{up to
multiples of $\omega_l$}$$
where $j=1,\ldots,l-1$.

On the other hand we have $\mathbb{B}_i([y^{i_1},\ldots,y^{i_q}])=
[x^{i_1n/l},\ldots,x^{i_qn/l}]$ which leads to:
$gcd(n,i_sn/l)=k_sn/l$, $s=1,\ldots,q$, $i_sn/l=k_su_sn/l$ and thus
$$r_{sj}^n\equiv ju_s\mod{n/(k_sn/l)=l/k_s}\;\;,\;\;j=1,\ldots,n-1$$
where $1\leq r_{sj}^n\leq n/(k_sn/l)=l/k_s$. By equation (\ref{n_j}) we obtain
for the multiplicities of the representations:
\begin{eqnarray*}
-n_j^n&\equiv& \big(\sum_{s=1}^qk_sr_{sj}^nn/l\big)/n\;,\;\text{up to
multiples of $\omega_n$} \\
&=&\big(\sum_{s=1}^qk_sr_{sj}^l\big)/l
\end{eqnarray*}
where $j=1,\ldots,n-1$.

This yields $n_j^l\equiv n_i^n$, up to multiples of the regular
representation, for all $j\equiv i\mod{l}$, $j=1,\ldots,l-1$,
$i=1,\ldots,n-1$. But this is exactly the
property of an induced representation, i.e.,
$\varphi(\phi_{\mathbb{B}_i\alpha})\equiv
Ind^G_H\varphi(\phi_{\alpha})$, up to multiples of $\omega_G$ and $\rho_0$.
Thus we have proven the assertion
$\theta\circ\mathbb{B}_i([y^{i_1},\ldots,y^{i_q}])=\widetilde{Ind}^G_H\circ
\theta([y^{i_1},\ldots,y^{i_q}])$.
\end{proof}
By the notation
$\bar{a}$ for an element $a\in\widetilde{R_{\mathbb{C}}G}$ we mean the
following. Let $\xi$ be a representative of the class $a$ in
$R_{\mathbb{C}}G$ then the complex conjugate $\bar{\xi}$ is a
representative of $\bar{a}$.
\begin{prop}\label{Lefschetz}
Let $\alpha\in\mathbb{B}_{C_n}$. Then the following holds.
$$\theta(-\alpha)=\overline{\theta(\alpha)}$$
\end{prop}
\begin{proof}
In the following we will use the notation of equation (\ref{n_j}). Let $\alpha=
[x^{i_1},\ldots,x^{i_q}]$, then $-\alpha=[x^{-i_1},\ldots,x^{-i_q}]$ thus 
$r^{\alpha}_{sj}\equiv jl_s \mod{n/k_s}$ and $r^{-\alpha}_{sj}\equiv 
j(-l_s)\equiv (n-j)l_s \mod{n/k_s}$. Insert this in equation (\ref{n_j}) and we
obtain
$$n^{\alpha}_j=n^{-\alpha}_{n-j}$$
which yields the proof.
\end{proof}
From proposition \ref{Lefschetz} we deduce that
$-\theta(\alpha)=\overline{\theta(\alpha)}$ and thus $\theta$ maps into 
$$\mathbb{A}_G:=\bigl\{a\in\widetilde{R_{
\mathbb{C}}G}|a+\bar{a}=0\bigr\}.$$
\begin{remark}
Proposition \ref{Lefschetz} follows also from the Lefschetz
Fixed-Point-Formula which states that the trace of $\varphi(\phi)+\overline{
\varphi(\phi)}$ equals two minus the number of fixed points.
\end{remark}
\begin{theo}\label{finiteindex}
The map $\theta:\mathbb{B}_{C_n}\rightarrow \widetilde{R_{\mathbb{C}}C_n}$ is 
injective and the index of the image of $\theta$ in $\mathbb{A}_{C_n}$ is 
finite.
\end{theo}
\begin{proof}
First we will prove the injectivity of $\theta$. Let $\phi$ be an
embedding of $C_n$ and $\alpha$ the 
corresponding element of $\phi$ in $\mathbb{B}_{C_n}$. \\
Assume that $\theta(\alpha)=0$. By lemma \ref{integralcyc} the representation
$\varphi(\phi)$ is integral and as the $C_n$-signature of Atiyah and
Singer \cite[p. 578]{AtSi} is given by
$sgn(\phi)=\varphi(\phi)-\overline{\varphi(\phi)}$ we have $sgn(\phi)=0$. By a
theorem of Edmonds and Ewing
\cite[Theorem 3.2]{EdEw} this implies that the fixed points of $\alpha$
are cancelling. (Note that this theorem doesn't say anything about the
other singular orbits.) We can now reduce $\alpha$ to a fixed point
free element of $\mathbb{B}_{C_n}$. \\
Let $H$ be a maximal proper subgroup of $C_n$. By
proposition \ref{Ressquare} we have
$\theta(\mathbb{B}res^G_H\alpha)=\widetilde{Res}^G_H\theta(\alpha)=0$. The same argument as
above tells us that the fixed points of $\mathbb{B}res^G_H\alpha$ cancel, i.e., 
the singular orbits of $\alpha$ which have $H$ as stabilizer cancel. We again
reduce $\alpha$ and iterate this procedure by restricting to smaller
and smaller subgroups of $C_n$. Finally we end up with a free action and
$\alpha=0$.\\
It is not difficult to see that $\mathbb{A}_{C_n}$ has rank $(n-1)/2$
for $n$ odd and $n/2-1$ for $n$ even. By examples \ref{exCmodd} and
\ref{exCmeven}, $\mathbb{B}_{C_n}$ has the same rank and as $\theta$ is
injective the index of $\theta(\mathbb{B}_{C_n})$ in
$\mathbb{A}_{C_n}$ is finite.
\end{proof}

\begin{prop}\label{indexp}
For $C_p$ the cyclic group of prime order $p$ the index of 
$\theta(\mathbb{B}_{C_p})$ in $\mathbb{A}_{C_p}$ is $h^-_p$, the 
first factor of the class number of $\mathbb{Q}(e^{2\pi i/p})$.
\end{prop}
\begin{proof}
The original proof can be found in Ewing's paper \cite[Theorem
3.2]{Ew1}. A proof which uses singular orbits and lies in the spirit
of this paper is given in \cite[Proposition 4]{Gr1}
\end{proof}
In view of proposition \ref{indexp} it is natural 
to ask if there is a relation between the finite index of theorem 
\ref{finiteindex} and the class number of $\mathbb{Q}(e^{2\pi i/n})$. This 
question is not settled for the moment but one can say the following. For $m$ 
an odd integer there exists a subgroup $M$ of $\mathbb{B}_{C_m}$ and a
subgroup $N$ of 
$\mathbb{A}_{C_m}$ such that the index of $M$ in $N$ is the class number of 
$\mathbb{Q}(e^{2\pi i/m})$.

\subsection{Finite Groups}
In this section we will prove similar properties for arbitrary finite groups 
as we did for cyclic groups in section \ref{cyclic}.

\begin{prop}
Let $\phi:G\hookrightarrow\G$ be an embedding. The corresponding
representation $\varphi(\phi):G\rightarrow U(g)$ can be computed by
$\sigma(\phi)=\bigl\{g;\hat{\gamma_1},\ldots,\hat{\gamma_q}\bigr\}_G$.
\end{prop}
\begin{proof}
Let $\{G_i\}_{i\in I}$ be the family of the maximal cyclic subgroups of
$G$. If a representation of $G$ is known when restricted to all the
maximal cyclic subgroups then the representation itself is known. 

We can restrict the extended data to any $G_i$ and we
obtain a data $Res^G_{G_i}\sigma(\phi)=\bigl\{
g;\hat{\beta_1},\ldots,\hat{\beta_n}\bigr\}_{G_i}$ which belongs to the
restricted embedding $\phi_i=\phi|_{G_i}:G_i\hookrightarrow
G\hookrightarrow\G$. 

By lemma \ref{repcyclic} we can compute $\varphi(\phi_i)$ for every
$G_i$ and thus the representation $\varphi(\phi)$ is given by the
extended data.
\end{proof}
In the following $\rho_0^G$ will always denote the one dimensional
trivial representation of $G$, usually we will just write $\rho_0$
when there is no confusion possible. 
\begin{example}\label{emptyfinite}
Let $G$ be a finite group and $\sigma(\phi)=\bigl\{g;\varnothing\bigr\}_G$ an 
extended data for a 
free action with Riemann-Hurwitz equation (\ref{RiHu}) $g=1+|G|(h-1)$. 
Restricting this extended data to any cyclic subgroup $H$ gives
$Res^G_H(\sigma(\phi))=\bigl\{g;\varnothing\bigr\}_H$ and Riemann-Hurwitz 
equation 
$g=1+|H|(\tilde{h}-1)$. Thus we have by example \ref{freesing} that 
$\varphi(\phi|_H)=\rho_0+$ \linebreak
$\omega_H (h-1)|G|/|H|$. We can now 
reconstruct the representation $\varphi(\phi)$ and obtain
$$\varphi(\phi)=\rho_0+ (h-1)\cdot\omega_G.$$
\end{example}
\begin{prop}\label{zerofinite}
The representation $\varphi(\phi)$ associated to $\sigma(\phi)=
\bigl\{g;$ $\hat{\gamma},\hat{\gamma}^{-1}\bigr\}_G$ for any finite 
group $G$ and any element $\gamma\in G$ equals
$$\xi=h\cdot \omega_G+\rho_0-Ind^G_{<\gamma>}\rho_0.$$
\end{prop}
\begin{proof}
Let $\nu$ denote the order of the element $\gamma$. Then we can write the 
Riemann-Hurwitz equation in the form $g=|G|h+1-|G|/\nu$. 
For $H$ a cyclic subgroup of $G$ let $t_i\in G$ $i=1,\ldots,|G|/|H|$ be the 
left coset 
representatives; then we can write $\beta_i=t_i\gamma^{r_i}t_i^{-1}$,
$r_i=\text{min}
\{s\in\mathbb{N}|t_i\gamma^s t_i^{-1}\in H\}$. We can renumber the 
$\beta_i$'s such that $\beta_i\neq 1$ for 
$i=1,\ldots,t$ and $\beta_i=1$ for $i=t+1,\ldots,|G|/|H|$.
Restricting $\phi$ to $H$ gives rise to the following extended
data and representation (see example \ref{zerosing}).
\begin{eqnarray*}
Res^G_H\sigma(\phi) & = & \bigl\{g;\hat{\beta_1},\hat{\beta_1}^{-1},\ldots,
\hat{\beta_t},\hat{\beta_t}^{-1}\bigr\}_H \\
\varphi(\phi|_H) 
& = &
\bigl((h-1)|G|/|H|+t\bigr)\cdot\omega_H+\rho_0-\sum^t_{i=1}Ind^H_{<\beta_i>}
\rho_0\\
& = & |G|/|H|\cdot h\cdot\omega_H+\rho_0-\sum^{|G|/|H|}_{i=1}
Ind^H_{<\beta_i>}\rho_0 \\
& = & Res^G_H\xi
\end{eqnarray*}
As $\varphi(\phi)$ and $\xi$ coincide on the cyclic subgroups they are
the same representation.
\end{proof}
We can define the following subgroups of $R_{\mathbb{C}}G$ 
\begin{eqnarray*}
D_G& := & \bl Ind^G_H\rho_0\,|\,H=G\;\text{and}\;H<G\;,\;H\;\text{cyclic}\br\\
E_G & := & \bl Ind^G_H\rho_0\,|\,H\leq G \br .
\end{eqnarray*}
By example \ref{emptyfinite} and proposition \ref{zerofinite},
$D_G$ is generated by representations $\varphi(\phi)$ where $\phi$
runs through free actions and actions which consist of cancelling
pairs. Notice that $D_G=E_G$ if $G$ is cyclic. But in general $D_G$ is 
strictly contained in $E_G$. In addition lemma \ref{integralcyc} is not true 
for arbitrary groups; thus we have
$$D_G<E_G<R_{\mathbb{Q}}G.$$
In the same way as for cyclic groups, example \ref{emptyfinite} shows
that there is a well defined map 
$$\Lambda_G\rightarrow R_{\mathbb{C}}G/\bigl<\rho_0,\omega_G\bigr>$$
induced by $\varphi$ and by proposition \ref{zerofinite} this map extends to
$$\tilde{\theta}:\mathbb{B}_G\rightarrow R_{\mathbb{C}}G/D_G.$$
In order to have commutative squares as in propositions
\ref{Ressquare} and \ref{Indsquare} we need another map
$$\theta:\mathbb{B}_G\rightarrow
R_{\mathbb{C}}G/E_G=:\widetilde{R_{\mathbb{C}}G}.$$ 
The map $\theta$ factors through $\tilde{\theta}$ and for $G$ a cyclic
group they are the same as we have in this case $D_G=E_G$.
For a subgroup $H$ of $G$ the restriction map
$\widetilde{Res}^G_H:\widetilde{R_{\mathbb{C}}G}\rightarrow \widetilde{R_{\mathbb{C}}H}$
and induction map $\widetilde{Ind}^G_H:\widetilde{R_{\mathbb{C}}H}\rightarrow
\widetilde{R_{\mathbb{C}}G}$ are again
well defined as $Res^G_HE_G\subset E_H$ and $Ind^G_HE_H\subset E_G$. 
\begin{prop}\label{GRessquare}
For a finite group $G$ and a subgroup $H<G$ the following diagram commutes.
\begin{eqnarray*}
\mathbb{B}_H & \stackrel{\mathbb{B}res^G_H}{\leftarrow} & \mathbb{B}_G \\
\downarrow & & \downarrow \\
\widetilde{R_{\mathbb{C}}H} & \stackrel{\widetilde{Res}^G_H}{\leftarrow} & \widetilde{R_{\mathbb{C}}G}
\end{eqnarray*}
\end{prop}
\begin{proof}
The proof is analogous to the one for proposition \ref{Ressquare}.
\end{proof}
\begin{lemma}\label{thetafs}
Let $G$ be any finite group with subgroup $H$ and $f_s:H\rightarrow
sHs^{-1}$ conjugation by $s\in G$. Then
there are induced maps $\mathbb{B}_{f_s}:\mathbb{B}_H\rightarrow
\mathbb{B}_{sHs^{-1}}$, $f_{s*}:R_{\mathbb{C}}H\rightarrow
R_{\mathbb{C}}sHs^{-1}$ 
and $\tilde{f}_{s*}:\widetilde{R_{\mathbb{C}}H}\rightarrow
\widetilde{R_{\mathbb{C}}sHs^{-1}}$. The map $f_{s*}$ is
given by $f_{s*}(\rho)(a)=\rho(s^{-1}as)$,
$\rho\in R_{\mathbb{C}}H$, $a\in sHs^{-1}$, and the map
$\tilde{f}_{s*}$ is induced by $f_{s*}$ on the quotient
$\widetilde{R_{\mathbb{C}}H}$. 
These maps satisfy the following equations:
\begin{eqnarray*}
\varphi(\phi_{\mathbb{B}_{f_s}\alpha})&\equiv&f_{s*}(\varphi(\phi_{\alpha}))
\;,\;\text{up to multiples of $\omega_{sHs^{-1}}$}\\
\theta\circ\mathbb{B}_{f_s}(\alpha)&=&\tilde{f}_{s*}\circ\theta(\alpha)
\end{eqnarray*}
for $\alpha\in\mathbb{B}_H$, $\mathbb{B}_{f_s}(\alpha)\in\mathbb{B}_{sHs^{-1}}$  
with corresponding embeddings $\phi_{\alpha}$ and $\phi_{\mathbb{B}_{f_s} 
\alpha}$ respectively.
\end{lemma}
\begin{proof}
The second equation follows directly from the first one. To prove the
first one just note that:
$$\varphi(\phi_{\mathbb{B}_{f_s}\alpha})\equiv\varphi(f_{s*}(\phi_{\alpha}))
=f_{s*}(\varphi(\phi_{\alpha}))\;,\;\text{up to multiples of
$\omega_{sHs^{-1}}$},$$
where of course $f_{s*}(\phi_{\alpha})(a)=\phi_{\alpha}(s^{-1}as)$.
\end{proof}
\begin{prop}\label{GIndsquare}
For a finite group $G$ and a subgroup $H<G$ the following diagram commutes.
\begin{eqnarray*}
\mathbb{B}_H & \stackrel{\mathbb{B}_i}{\rightarrow} & \mathbb{B}_G \\
\downarrow &  & \downarrow \\
\widetilde{R_{\mathbb{C}}H} & \stackrel{\widetilde{Ind}^G_H}{\rightarrow} & \widetilde{R_{
\mathbb{C}}G}
\end{eqnarray*}
\end{prop}
\begin{proof}
To prove the statement of the proposition it is enough to prove the
following equation:
\begin{eqnarray}\label{GInd}
\varphi(\phi_{\mathbb{B}_i\alpha})\equiv
Ind^G_H\varphi(\phi_{\alpha})\;,\;
\text{up to multiples of $\omega_G$ and $\rho_0^G$}
\end{eqnarray}
for some elements $\mathbb{B}_i\alpha\in\mathbb{B}_G$,
$\alpha\in\mathbb{B}_H$ with
corresponding embeddings $\phi_{\mathbb{B}_i\alpha}$ and
$\phi_{\alpha}$ respectively.
To do this, we will restrict both sides of equivalence (\ref{GInd}) to an
arbitrary cyclic subgroup $K$ of $G$ and then show that they 
coincide, up to multiples of $\omega_K$ and $\rho_0^K$. This will then
imply that the original equivalence (\ref{GInd}) is correct.

Let $K$ be a cyclic subgroup of $G$ and $S$ a set of representatives
for the $(H,K)$ double cosets of $G$. For $s\in S$, let
$H_s=sHs^{-1}\cap K$, which is a subgroup of $K$. The inclusions
will be denoted by $i:H\rightarrow G$ and $j:H_s\rightarrow K$. $f_s$
will be conjugation by $s$, i.e., $f_s(g)=sgs^{-1}$. 

In the following we will write $\varphi\circ F(\alpha)$ instead of
$\varphi(\phi_{F(\alpha}))$ for any map $F$. All the equivalences will be
up to multiples of $\omega_K$ and $\rho_0^K$. 

First we will apply $Res^G_K$ to the right hand side of equivalence
(\ref{GInd}). Using proposition \ref{GRessquare}, lemma \ref{thetafs}
and the well known double coset formula for representations we obtain
\begin{eqnarray*}
Res^G_K\circ Ind^G_H\circ\varphi &=& \bigoplus_{s\in
S}Ind^K_{H_s}\circ Res^{sHs^{-1}}_{H_s}\circ f_{s*}\circ\varphi \\
& \equiv &\bigoplus_{s\in S}Ind^K_{H_s}\circ\varphi\circ
\mathbb{B}res^{sHs^{-1}}_{H_s}\circ\mathbb{B}_{f_s}.
\end{eqnarray*}
Now we apply $Res^G_K$ to the left hand side of equivalence
(\ref{GInd}). Using propositions \ref{GRessquare} and
\ref{doublecoset} we obtain
\begin{eqnarray*}
Res^G_K\circ\varphi\circ\mathbb{B}_i &\equiv &\varphi\circ
\mathbb{B}res^G_K\circ\mathbb{B}_i\equiv\varphi\circ\bigoplus_{s\in S}\mathbb{B}_j\circ
\mathbb{B}res^{sHs^{-1}}_{H_s}\circ\mathbb{B}_{f_s}\\
&\equiv& \bigoplus_{s\in S}\varphi\circ\mathbb{B}_j\circ
\mathbb{B}res^{sHs^{-1}}_{H_s}\circ\mathbb{B}_{f_s}.
\end{eqnarray*}
By the proof of proposition \ref{Indsquare} we know that
$\varphi\circ\mathbb{B}_j\equiv Ind^K_{H_s}\circ\varphi$, up to
multiples of $\omega_K$, for cyclic groups $K$. Thus both sides of
equivalence (\ref{GInd}) coincide up to multiples of $\omega_G$ and
$\rho_0^G$ and the proposition is true. 
\end{proof}
Recall from corollary \ref{ZZ/2Z} that $\mathbb{B}_G$ consists only of
copies of $\mathbb{Z}$ and $\mathbb{Z}/2\mathbb{Z}$.
\begin{theo}\label{Zinj}
The map $\theta$ is injective on the copies of $\mathbb{Z}$ in $\mathbb{B}_G$.
\end{theo}
\begin{proof}
Let $\bigl[\hat{\gamma_1},\ldots,\hat{\gamma_q}\bigr]$ be a reduced element
of infinite order in $\mathbb{B}_G$, i.e., $\gamma_i$ is not conjugate
to $\gamma_j^{-1}$, $i\neq j$, $i,j=1,\ldots,q$ and there exists at
least one $\gamma_{i_0}$, $1\leq i_0\leq q$, which is not conjugate to its
inverse. 
There is also an $s=1,\ldots,q$ such that we can renumber the
$\gamma_i$'s in the following way: $i_0=1$,
$\bigl<\gamma_i\bigr>$ is conjugate to $\bigl<\gamma_1\bigr>$, $i\leq
s\leq q$, and
$\bigl<\gamma_j\bigr>$ is not conjugate to $\bigl<\gamma_1\bigr>$,
$s<j\leq q$. Then we have that $\gamma_i$ is not conjugate to
$\gamma_j^{-1}$, $1\leq i,j\leq s$.

 Notice that $\gamma_1$
can be chosen such that $\bigl<\gamma_1\bigr>$ is maximal with respect to the other
cyclic subgroups $\bigl<\gamma_i\bigr>$, $i=s+1,\ldots,q$, up to conjugation.

The singular orbit data
$\alpha=\mathbb{B}res^G_{<\gamma_1>}\bigl[\hat{\gamma_1},\ldots,\hat{\gamma_q}\bigr]$ has
then fixed
points which don't cancel and thus $\alpha\neq0$. (Notice that the
other singular orbits 
of $\alpha$ could cancel.) By theorem \ref{finiteindex} we know that
$\theta(\alpha)\neq0$ and by proposition \ref{GRessquare} we deduce
$$\widetilde{Res}^G_{<\gamma_1>}\theta(\bigl[\hat{\gamma_1},\ldots,\hat{\gamma_q}\bigr])=
\theta(\alpha)\neq0$$
and thus
$\theta(\bigl[\hat{\gamma_1},\ldots,\hat{\gamma_q}\bigr])\neq0$
\end{proof}
\begin{Kor}\label{Gabelian}
Let $G$ be an Abelian group which doesn't contain a copy of
$\mathbb{Z}/2\mathbb{Z}\times\mathbb{Z}/2\mathbb{Z}$. The map
$\theta$ is then injective. 
\end{Kor}
\begin{proof}
By corollary \ref{BGtorsion}, $\mathbb{B}_G$ is torsion free and 
by theorem \ref{Zinj}, $\theta$ is injective on the copies of
$\mathbb{Z}$.
\end{proof}
\begin{prop}\label{AG}
For $G$ a finite group, $\theta:\mathbb{B}_G\rightarrow
\widetilde{R_{\mathbb{C}}G}$ maps into 
$$\mathbb{A}_G:=\bigl\{a\in\widetilde{R_{
\mathbb{C}}G}|a+\bar{a}=0\bigr\},$$
i.e., $\theta(-\alpha)=\overline{\theta(\alpha)}$. 
\end{prop}
\begin{proof}
Let $\alpha$ be an element of $\mathbb{B}_G$ with corresponding embedding
$\phi_{\alpha}$ and unitary representation
$\varphi(\phi_{\alpha})$. Changing the orientation of the surface
is equivalent with taking the transposed symplectic structure in the
first homology group of the surface. On the other hand transposition
in $Sp_{2n}(\mathbb{R})$ induces complex conjugation in its unitary
subgroup $U(n)\subset Sp_{2n}(\mathbb{R})$ .

Thus to the element $-\alpha$ corresponds the unitary representation
$\overline{\varphi(\phi_{\alpha})}$ and we obtain 
$-\theta(\alpha)=\theta(-\alpha)=\overline{\theta(\alpha)}$.
\end{proof}
\begin{Kor}\label{repconj}
For an element $\alpha$ of order two in $\mathbb{B}_G$ we have
$$\theta(\alpha)=\overline{\theta(\alpha)}$$
and for an element $\beta$ of infinite order 
$$\theta(\beta)\neq\overline{\theta(\beta)}.$$
\end{Kor}
\begin{proof}
As $\alpha=-\alpha$ we have by proposition \ref{AG},
$\theta(\alpha)=\theta(-\alpha)=\overline{\theta(\alpha)}$. 

In the second case it follows that $\beta\neq-\beta$ and thus by theorem
\ref{Zinj} and proposition \ref{AG},
$\theta(\beta)\neq\theta(-\beta)=\overline{\theta(\beta)}$.  
\end{proof}

\begin{remark}\label{S_3}
For $G=S_3$, the symmetric group on three letters, we have 
$D_G=R_{\mathbb{Q}}G=R_{\mathbb{C}}G$. Thus $\theta$ is the zero map. 

Note that the $G$-signature defined by Atiyah and Singer
$sgn(\phi)=\varphi(\phi)-\overline{\varphi(\phi)}$  is also zero. But with
the approach we chose, it is now possible to modify $\theta$ such
that it becomes injective again.

For $S_3$ let $\chi_0$ be the trivial representation, $\chi_1$ the one
dimensional non-trivial representation and $\chi_2$ the
two dimensional irreducible representation. Recall the notation from
example \ref{exS_3} and let 
$$D'_{S_3}=\bigl<Ind^{S_3}_{<a>}\rho_0\,,\,Ind^{S_3}_{<1>}\rho_0\,,\,Ind^{S_3}_{S_3}\rho_0\bigr>
=\bigl<\chi_0+\chi_1\,,\,\chi_0+\chi_1+2\chi_2\,,\,\chi_0\bigr>.$$
Then we have that the class of $\chi_2$ generates
$R_{\mathbb{C}}S_3/D'_{S_3}\cong C_2$. We know
that $\mathbb{B}_{S_3}\cong\bigl<\bigl[\hat{a}\bigr]\bigr>\cong C_2$; thus it
remains to show that $\theta':\mathbb{B}_{S_3}\rightarrow
R_{\mathbb{C}}S_3/D'_{S_3}$ is non zero on $\bigl[\hat{a}\bigr]$.

Let $\sigma(\phi)=\bigl\{3;\hat{a}\bigr\}_{S_3}$ then we find that
$\varphi(\phi)=\chi_0+\chi_2$ and thus
$\theta'$ maps $\bigl[\hat{a}\bigr]$ to the class of $\chi_2$ and
$\theta'$ is an isomorphism.

Note that $D'_{S_3}$ consists of the minimal relations in order to make
$\theta'$ a group homomorphism. It consists of the representations
coming from the free actions and the cancelling pair
$\bigl[\hat{a},\hat{a}\bigr]_{S_3}$.

In this situation the maps
$\mathbb{B}_i:\mathbb{B}_{C_3}\rightarrow\mathbb{B}_{S_3}$,
$\widetilde{Ind}_{C_3}^{S_3}:\widetilde{R_{\mathbb{C}}C_3}\rightarrow
R_{\mathbb{C}}S_3/D'_{S_3}$ and $\theta'$ still commute. This follows
from the fact that  
$$Ind_{C_3}^{S_3}(\varphi(\phi_{[a,a,a]_{C_3}}))\equiv
Ind_{C_3}^{S_3}\rho_1=\chi_2\equiv\varphi(\phi_{[\hat{a}]_{S_3}})
\equiv
\varphi(\phi_{\mathbb{B}_i([a,a,a]_{C_3}}))$$ 
where the equivalence is taken modulo $D'_{S_3}$. $\rho_1$ denotes a
one dimensional faithful representation of $C_3$. 
\end{remark}

\section{The Case $\mathbb{Z}/p\mathbb{Z}\times\mathbb{Z}/p\mathbb{Z}$}
\label{case}
In this section we want to have a closer look at the case
$G=\mathbb{Z}/p\mathbb{Z}\times\mathbb{Z}/p\mathbb{Z}$, $p$ an odd
prime. We know by example \ref{CpCp} that $\mathbb{B}_G$ is isomorphic to
$\mathbb{Z}^{\frac{p^2-1}{2}}$ and thus by theorem \ref{Zinj} and
proposition \ref{AG},
$\theta:\mathbb{B}_G\rightarrow\mathbb{A}_G$ is injective. The main
purpose of this section is to prove the following theorem.
\begin{theo}\label{CpCpindex}
The index $\Delta$ of $\theta(\mathbb{B}_G)$ in $\mathbb{A}_G$ is
$(h^-_p)^{p+1}\cdot p^i$, where $i$ is an integer in the range $1-p+k\leq
i\leq k+(p-1)^2/2$ for some $k\in\mathbb{N}$; $h^-_p$ denotes the first factor of the class number
of $\mathbb{Q}(e^{2\pi i/p})$; i.e., $\Delta=\bigl[\mathbb{A}_G:\theta(\mathbb{B}_G)\bigr]=(h^-_p)^{p+1}\cdot p^i$.
\end{theo}
Before we can prove this theorem, we need some lemmata. Let $G_j$,
$j=0,\ldots,p$, be the distinct cyclic subgroups of order $p$ in
$G$. In the sequel $\mathbb{B}_i$ will denote a map
$\mathbb{B}_{G_j}\rightarrow \mathbb{B}_G$ coming from an inclusion
$i:G_j\rightarrow G$. We can consider the following commutative diagram.
\begin{equation*}
\CD
\bigoplus\mathbb{B}_{G_j} @>\oplus\mathbb{B}_i>> \mathbb{B}_G
@>\oplus \mathbb{B}res^G_{G_j}>>  \bigoplus\mathbb{B}_{G_j} \\
@V\oplus\theta VV   @V\theta VV @V\oplus\theta VV \\
\bigoplus\mathbb{A}_{G_j} @>\oplus \widetilde{Ind}^G_{G_j}>> \mathbb{A}_G
@>\oplus \widetilde{Res}^G_{G_j}>>  \bigoplus\mathbb{A}_{G_j} 
\endCD
\end{equation*}
{\it Notation}: 
In the sequel if we want to talk about the index of, for example, 
$\oplus \mathbb{B}res^G_{G_j}\circ\oplus\mathbb{B}_i(\bigoplus\mathbb{B}_{G_j})$
in $\bigoplus\mathbb{B}_{G_j}$ we will omit the maps and just
write the index of $\bigoplus\mathbb{B}_{G_j}$ in
$\bigoplus\mathbb{B}_{G_j}$.\\[1mm] 
By theorem \ref{finiteindex} and
proposition \ref{indexp} the map $\oplus\theta$ is injective and the
index of $\bigoplus\mathbb{B}_{G_j}$ in $\bigoplus\mathbb{A}_{G_j}$ is
$(h^-_p)^{p+1}$. An easy computation shows that the map $\oplus
\mathbb{B}res^G_{G_j}\circ\oplus\mathbb{B}_i$ is just multiplication by $p$,
thus $\oplus\mathbb{B}_i$ is injective and the index of
$\bigoplus\mathbb{B}_{G_j}$ in $\bigoplus\mathbb{B}_{G_j}$ is
$p^{(p^2-1)/2}$. (By example \ref{exCmodd},
$\mathbb{B}_{G_j}\cong\mathbb{Z}^{(p-1)/2}$.) On the other hand we
have also that the map $\oplus\mathbb{B}_i\circ\oplus \mathbb{B}res^G_{G_j}$ is
multiplication by $p$ and thus $\oplus \mathbb{B}res^G_{G_j}$ is injective and
the index of $\mathbb{B}_G$ in $\mathbb{B}_G$ is also $p^{(p^2-1)/2}$.
\begin{lemma}\label{Bgjindex}
The index of $\bigoplus\mathbb{B}_{G_j}$ in $\mathbb{B}_G$ is $p^{p-1}$.
\end{lemma}
\begin{proof}
Recall the notations of examples \ref{exCmodd} and \ref{CpCp} and let
$G_j=\bigl<xy^j\bigr>$, $j=0,\ldots,p-1$ and $G_p=\bigl<y\bigr>$. 
A basis for $\mathbb{B}_{G_j}$ is given by:
\begin{align*}
[x,x^k,x^{p-k-1}]_{G_0}&\;,\;k=1,\ldots,(p-1)/2\\
[y,y^l,y^{p-l-1}]_{G_p}&\;,\;l=1,\ldots,(p-1)/2\\
\bigl[xy^j,(xy^j)^r,(xy^j)^{p-r-1}\bigr]_{G_j}&\;,\;r=1,\ldots,(p-1)/2\;,\;j=1,\ldots,p-1.
\end{align*}
For $\mathbb{B}_G$ we have the following basis:
\begin{align*}
[x,x^k,x^{p-k-1}]_G & =\mathbb{B}_i([x,x^k,x^{p-k-1}]_{G_0})\;,\;
k=1,\ldots,(p-1)/2 \\
[y,y^l,y^{p-l-1}]_G & =\mathbb{B}_i([y,y^l,y^{p-l-1}]_{G_p})\;,\;
l=1,\ldots,(p-1)/2 \\
[x^j,y^i,x^{p-j}y^{p-i}]_G &\;,\; j=1,\ldots,(p-1)/2\;,\;i=1,\ldots,p-1.
\end{align*}
Thus we see that the basis elements of $\mathbb{B}_G$,
$\bigl[x,x^k,x^{p-k-1}\bigr]_G$, $k=1,\ldots,(p-1)/2$, and
$\bigl[y,y^l,y^{p-l-1}\bigr]_G$, $l=1,\ldots,(p-1)/2$, are induced form
$\mathbb{B}_{G_0}$ and $\mathbb{B}_{G_p}$. Consequently they are zero
in the quotient $Q=\mathbb{B}_G/\bigoplus\mathbb{B}_{G_j}$. It remains
to find the relations between the other basis elements of
$\mathbb{B}_G$ in $Q$.

We can write the images of the basis elements of $\mathbb{B}_{G_j}$,
$j=1,\ldots,p-1$, as a linear combination in the basis elements of
$\mathbb{B}_G$.
\begin{gather*}
\mathbb{B}_i(\bigl[xy^j,(xy^j)^r,(xy^j)^{p-r-1}\bigr]_{G_j})=
\ominus\bigl[x,y^j,(xy^j)^{p-1}\bigr]\ominus\bigl[x^r,y^{rj},(xy^j)^{p-r}\bigr]
\\
\oplus\bigl[x^{r+1},y^{j(r+1)},(xy^j)^{p-r-1}\bigr]
\oplus\bigl[x,x^r,x^{p-r-1}\bigr]_G\oplus\bigl[y^j,y^{jr},y^{p-jr-j}\bigr]_G
\end{gather*}
(Note that $\bigl[y^j,y^{jr},y^{p-jr-j}\bigr]_G$ can be written as a
linear combination in the $\bigl[y,y^l,y^{p-l-1}\bigr]_G$,
$l=1,\ldots,(p-1)/2$.) If we
consider the above equation in the quotient $Q$, we obtain the
remaining relations between the basis elements of $\mathbb{B}_G$ in $Q$:
$$\bigl[x,y^j,(xy^j)^{p-1}\bigr]_G\oplus\bigl[x^r,y^{rj},(xy^j)^{p-r}\bigr]_G
\equiv\bigl[x^{r+1},y^{j(r+1)},(xy^j)^{p-r-1}\bigr]_G$$
$r=1,\ldots,(p-1)/2$, $j=1,\ldots,p-1$.\\
From this we conclude that the quotient $Q$ is generated by the
elements $\bigl[x,y^j,(xy^j)^{p-1}\bigr]_G$, $j=1,\ldots,p-1$ and thus
the index is $p^{p-1}$.
\end{proof}
From this lemma we deduce immediately that the index of $\mathbb{B}_G$
in $\bigoplus\mathbb{B}_{G_j}$ is $p^{(p-1)^2/2}$.
\begin{lemma}
The maps $\oplus \widetilde{Res}^G_{G_j}\circ\oplus
\widetilde{Ind}^G_{G_j}:\bigoplus\mathbb{A}_{G_j}\rightarrow\bigoplus\mathbb{A}_{G_j}$
and $\oplus \widetilde{Ind}^G_{G_j}\circ\oplus
\widetilde{Res}^G_{G_j}:\mathbb{A}_G\rightarrow\mathbb{A}_G$ are just
multiplication by $p$.
\end{lemma}
\begin{proof}
Take a class $a_j\in\mathbb{A}_{G_j}$ with representative
$\xi_j\in R_{\mathbb{C}}G_j$. The induced
representation $Ind^G_{G_j}(\xi_j)$ consists of $p$ copies of $\xi_j$
and $G/G_j$ operates by permutation on the $p$ copies. Thus $Res^G_{G_l}\circ
Ind^G_{G_j}(\xi_j)$ is $p\cdot a_j$ in $\mathbb{A}_{G_l}$ whenever $l=j$
and zero when $l\neq j$. This proves the assertion for the first
map.

For the second map take a class $a\in\mathbb{A}_G$ with representative
$\xi\in R_{\mathbb{C}}G$. Then restrict $\oplus
Ind^G_{G_j}\circ\oplus Res^G_{G_j}(\xi)$ to a subgroup $G_l$ of
$G$ and we obtain:
\begin{eqnarray*}
Res^G_{G_l}\circ\oplus Ind^G_{G_j}\circ\oplus Res^G_{G_j}(\xi) &=&
Res^G_{G_l}\circ\oplus_j(Ind^G_{G_j}\circ Res^G_{G_j})(\xi)= \\
\oplus_j Res^G_{G_l}\circ Ind^G_{G_j}\circ Res^G_{G_j}(\xi) &=&
p\cdot Res^G_{G_l}(\xi)\oplus\text{multiples of $\omega_{G_l}$}.
\end{eqnarray*}
The map $\oplus_l Res^G_{G_l}$ is injective on $R_{\mathbb{C}}G$ 
thus the map $\oplus
Ind^G_{G_j}\circ\oplus Res^G_{G_j}$ acts by multiplication by $p$, up
to multiples of $\omega_G$ on $R_{\mathbb{C}}G$ and
consequently the map $\oplus \widetilde{Ind}^G_{G_j}\circ\oplus \widetilde{Res}^G_{G_j}$ acts
by multiplication by $p$ on $\mathbb{A}_G$.
\end{proof}
An immediate consequence of this lemma is that the map $\oplus
\widetilde{Ind}^G_{G_j}$ is injective and that the index
of $\bigoplus\mathbb{A}_{G_j}$ in $\bigoplus\mathbb{A}_{G_j}$ is
$p^{(p^2-1)/2}$. (Note that
$\mathbb{A}_{G_j}\cong\mathbb{Z}^{(p-1)/2}$.) Furthermore the kernel
of the map $\oplus \widetilde{Res}^G_{G_j}:\mathbb{A}_G\rightarrow\bigoplus
\mathbb{A}_{G_j}$ is isomorphic to $(\mathbb{Z}/p\mathbb{Z})^k$ for
some $k\in\mathbb{N}$. 
\begin{lemma}\label{isoAG}
The group $\mathbb{A}_G$ is isomorphic to $\mathbb{Z}^{\frac{p^2-1}{2}}
\oplus(\mathbb{Z}/p\mathbb{Z})^k$ for some $k\in\mathbb{N}$.
\end{lemma}
\begin{proof}
We know that $\bigoplus\mathbb{A}_{G_j}$ is isomorphic to
$\mathbb{Z}^{(p^2-1)/2}$ and that it injects into $\mathbb{A}_G$. This
proves that $\mathbb{Z}^{(p^2-1)/2}$ is a direct summand of
$\mathbb{A}_G$. On the other hand the kernel of $\oplus \widetilde{Res}^G_{G_j}:
\mathbb{A}_G\rightarrow\bigoplus\mathbb{A}_{G_j}\cong\mathbb{Z}^{(p^2-1)/2}$
is isomorphic to $(\mathbb{Z}/p\mathbb{Z})^k$ for some
$k\in\mathbb{N}$ and thus the assertion holds.
\end{proof}
{\it Proof of theorem \ref{CpCpindex}.} 
By example \ref{CpCp} the group $\mathbb{B}_G$ is isomorphic to
$\mathbb{Z}^{\frac{p^2-1}{2}}$ and thus by theorem \ref{Zinj} the map
$\theta:\mathbb{B}_G\rightarrow\mathbb{A}_G$ is injective. Lemma
\ref{isoAG} tells us now that the index $\Delta$ is finite.

There are two possibilities to compute the index
$\bigl[\mathbb{A}_G:\bigoplus\mathbb{B}_{G_j}\bigr]$. For the first
possibility we use lemma \ref{Bgjindex}.
$$\bigl[\mathbb{A}_G:\bigoplus\mathbb{B}_{G_j}\bigr]=\bigl[\mathbb{A}_G:
\mathbb{B}_G\bigr]\cdot\bigl[\mathbb{B}_G:\bigoplus\mathbb{B}_{G_j}\bigr]=
\Delta\cdot p^{p-1}$$
From lemma \ref{isoAG} we can also deduce that the index of
$\bigoplus\mathbb{A}_{G_j}$ in $\mathbb{A}_G$ is some prime power
$p^i$ with $i=0,\ldots,k+(p^2-1)/2$, thus for the second possibility we
get
\begin{eqnarray*}
\bigl[\mathbb{A}_G:\bigoplus\mathbb{B}_{G_j}\bigr]&=&\bigl[\mathbb{A}_G:
\bigoplus\mathbb{A}_{G_j}\bigr]\cdot\bigl[\bigoplus\mathbb{A}_{G_j}:\bigoplus
\mathbb{B}_{G_j}\bigr]\\
&=&p^i\cdot(h^-_p)^{p+1}\;\;,\;\;i=0,\ldots,k+(p^2-1)/2.
\end{eqnarray*}
Collecting the results of the two equations we obtain
$\Delta=(h^-_p)^{p+1}\cdot p^i$ where $i$ ranges over the integers
$1-p+k\leq i\leq k+(p-1)^2/2$. \hfill$\square$
\begin{remark}\label{p=3}
In the case when the prime $p$ is $3$ the index $\bigl[\mathbb{A}_G:
\bigoplus\mathbb{A}_{G_j}\bigr]$ is also $p^{p-1}=9$ thus $\Delta=(h^-_3)^4=1$.
\end{remark}
In view of remark \ref{p=3} we can formulate a conjecture.
\begin{equation}\tag{Conjecture}
\Delta=(h^-_p)^{p+1}
\end{equation}
To prove this conjecture it would be enough to show that $\bigl[\mathbb{A}_G:
\bigoplus\mathbb{A}_{G_j}\bigr]$ $=p^{p-1}$. In general this index is
hard to compute, however for the prime $3$ we have seen that the
conjecture is true.

For regular primes $p$, i.e., $p$ doesn't divide its class
number $h^-_p$, we can give an equivalent statement to the conjecture. 
\begin{prop}\label{conjequiv}
Let $p$ be a regular prime. Statements 1. and 2. are equivalent.
\begin{enumerate}
\item 
\begin{enumerate}
\item
Let $a$ be an element of $\mathbb{A}_G$ such that its restriction to
every subgroup $G_j$ lies in the image of $\theta$, i.e.,
$\widetilde{Res}^G_{G_j}a\in\theta(\mathbb{B}_{G_j})$, $j=0,\ldots,p$. Then the
element $a$ itself lies already in $\theta(\mathbb{B}_G)$.
\item
Let $a_j$, $j=0,\ldots,p$, be elements in $\mathbb{A}_{G_j}$,
such that the sum of their induced images lies in
$\theta(\mathbb{B}_G)$, i.e., $\oplus
\widetilde{Ind}^G_{G_j}a_j\in\theta(\mathbb{B}_G)$. Then the representations
$a_j$ lie already in $\theta(\mathbb{B}_{G_j})$.
\end{enumerate}
\item
The conjecture $\Delta=(h^-_p)^{p+1}$ is true.
\end{enumerate}
\end{prop}
\begin{proof}
First we will prove 1.$\Rightarrow$2. .

Statement 1.(a) is equivalent to the fact that the induced map
$\mathbb{A}_G/\mathbb{B}_G$ $\rightarrow\bigoplus\mathbb{A}_{G_j}/\bigoplus
\mathbb{B}_{G_j}$ is injective. Statement 1.(b) is equivalent to the
fact that the induced map $\bigoplus\mathbb{A}_{G_j}/\bigoplus\mathbb{B}_{G_j}
\rightarrow\mathbb{A}_G/\mathbb{B}_G$ is injective. From this we
conclude $\Delta|(h^-_p)^{p+1}$ and $(h^-_p)^{p+1}|\Delta$ which in
turn proves the conjecture.

Next we will prove 2.$\Rightarrow$1. .

The map $\oplus \widetilde{Ind}^G_{G_j}\circ\oplus \widetilde{Res}^G_{G_j}$ is multiplication
by $p$. Thus the induced map
$\mathbb{A}_G/\mathbb{B}_G\rightarrow\mathbb{A}_G/\mathbb{B}_G$ is
injective as $p$ doesn't divide the class number $h^-_p$. From this we
conclude that the map $\mathbb{A}_G/\mathbb{B}_G\rightarrow\bigoplus
\mathbb{A}_{G_j}/\bigoplus\mathbb{B}_{G_j}$, which is induced by
$\oplus \widetilde{Res}^G_{G_j}$, is injective, which is equivalent to 1.(a).

The same argument using the map $\oplus \widetilde{Res}^G_{G_j}\circ\oplus
\widetilde{Ind}^G_{G_j}$ shows that 1.(b) is also true.
\end{proof}
\begin{remark}
Statement 1. of proposition \ref{conjequiv} is never true if the prime
$p$ is irregular.
\end{remark}

\section{$G$-Equivariant Cobordism}\label{cobordism}
In this section we will define a map $\chi$ from the $G$-equivariant
cobordism group of surface diffeomorphisms $\Omega_G$ to $\mathbb{B}_G$
and prove that this map is surjective and describe its kernel. We will
only talk about oriented cobordism, thus every diffeomorphism is
orientation preserving. First we need some definitions.
\begin{defin}
Let $M_1$ (resp. $M_2$) be a compact, oriented, connected Riemann
surface with smooth $G$-action $\kappa_1:G\rightarrow\text{\it
Diffeo}_+(M_1)$ (resp. $\kappa_2:G\rightarrow\text{\it Diffeo}_+(M_2)$) 
We say that $\kappa_1$ is $G$-equivariant cobordant to $\kappa_2$, written
$\kappa_1\sim\kappa_2$, if there exists a smooth, compact, oriented, connected
3-manifold $V$ and a smooth $G$-action $\Phi$ on $V$ such that
\begin{enumerate}
\item The boundary of $V$ is the disjoint union of $M_1$ and $-M_2$,
$\partial(V)=M_1\cup-M_2$. The notation $-M_2$ denotes $M_2$ with opposite
orientation. The orientations on $M_1$ and $-M_2$ coincide with the
one induced by $V$. 
\item $\Phi$ restricted to $\partial(V)$ agrees with $\kappa_1\cup\kappa_2$.
\end{enumerate}
We also say that $\kappa_1$ is zero $G$-equivariant cobordant, written
$\kappa_1\sim0$, if $\partial(V)=M_1$.
$\Omega_G$ will denote the set of $G$-equivariant cobordism classes
and a class will be denoted by $(\kappa,M)$.
\end{defin}
There is an addition in $\Omega_G$ given by the $G$-equivariant connected sum
defined in section \ref{geometry}. This addition is well defined as
it doesn't depend on the representative of the $G$-equivariant
cobordism classes, thus we have an Abelian monoid. In fact we have
even more, $\Omega_G$ is an Abelian group. Indeed, for any element 
$(\kappa,M)$, there exists also its inverse element $(\kappa,-M)$ in
$\Omega_G$. The zero element is the zero $G$-equivariant
cobordism class. 
\begin{lemma}\label{cobordef}
Let $\kappa_1$ and $\kappa_2$ be $G$-equivariant cobordant. Then they
have the same singular orbit data. 
\end{lemma}
\begin{proof}
There is a compact, connected smooth $3$-manifold $V$ together with a smooth
$G$-action $\Phi$ on $V$  such that $\partial(V)=M_1\cup-M_2$ and
$\Phi$ restricted to $\partial(V)$ agrees with
$\kappa_1\cup\kappa_2$. \\
Let $x\in M_1$ and $Gx$ be a singular orbit such that $\gamma\in G$
generates the stabilizer of $x$ and acts by rotation through
$2\pi/|\bigl<\gamma\bigr>|$.  
The action of $H=\bigl<\gamma\bigr>$ extends to $V$ by extending the
fixed point 
$x$ to a properly embedded arc of fixed points of $H$. The endpoint $y$ of
this arc has to lie in a boundary component, thus either in $M_1$ or in
$M_2$. If $y$ lies in $M_1$ we have again two possibilities.

Firstly $y$ lies in the same singular orbit $Gx$, but then there exists
$g\in G$ with $gx=y$ and thus $g$ maps the arc onto itself and
the diffeomorphism given by $g$ isn't orientation preserving anymore.
Secondly $y$ lies in another singular orbit $Gy$ and thus $\gamma$ acts by
rotation through $-2\pi/|\bigl<\gamma\bigr>|$ on $y$. The arc from $x$
to $y$ is mapped by $g$ to arcs from $gx$ to $gy$. Thus $Gx$ and $Gy$
are two singular orbits on $M_1$ which cancel.

If $y$ lies in $M_2$ then by the same arguments as above and the fact
that we have the opposite orientation on $M_2$ we deduce that $Gx$ and $Gy$
are the same singular orbits.

In this way we can show that up to cancelling pairs $\kappa_1$ and
$\kappa_2$ have the same singular orbits.
\end{proof}
By lemma \ref{cobordef} we can now define a homomorphism $\chi$ which
sends a $G$-equivariant cobordism class to its singular orbit data.
$$\chi:\Omega_G\rightarrow\mathbb{B}_G.$$
\begin{theo}\label{chi}
The map $\chi$ is surjective and the kernel consists of cobordism
classes of free $G$-actions.
\end{theo}
\begin{proof}
First we will prove the statement about the kernel. Let $\kappa$
denote a $G$-action on
a surface $M$ whose image in $\mathbb{B}_G$ is zero, i.e., the
singular orbits cancel. We will show that $\kappa$ is cobordant to a free
$G$-action. 

First we take the product cobordism
$V=M\times[0,1]$, where $G$ is extended over $V$ in the obvious
way. Next we modify $V$ on the top end $M\times\{1\}$ in a similar way
as we did to introduce relations in $W_G$ in section
\ref{geometry}. Suppose $Gx$ and $Gy$ are cancelling singular orbits
on $M\times\{1\}$ with stabilizer $\bigl<\gamma\bigr>$ at $x$ and $y$. Let $T$
be a set of
representatives for the $\bigl<\gamma\bigr>$ left cosets of $G$. Find discs
$D_1$ and $D_2$ around $x$ and $y$ respectively such that $D_j$ is fixed by
$\bigl<\gamma\bigr>$, $j=1,2$, and $\cup_{j=1,2}\cup_{t\in
T}\bigl\{tD_j\bigr\}$ are mutually disjoint. Then connect each pair of
discs $tD_1$, $tD_2$ by a solid tube $D^2\times[0,1]$ for every
element $t\in T$. Let $V_1$ denote the
3-manifold which results after connecting all the cancelling pairs of
singular orbits by solid tubes. The action of $G$ can be extended to $V_1$ by
rotating and permuting the solid tubes. The manifold $V_1$ provides
the $G$-equivariant cobordism between $\kappa$ and a free
$G$-action. 

To prove that $\chi$ is surjective just note that every singular
orbit data gives rise to a $G$-action with the appropriate singular
orbits and thus to the corresponding $G$-equivariant cobordism class.
\end{proof}
With the surjection $\chi$ the map $\theta$ becomes a
$G$-signature. It has common properties with the $G$-signature of
Atiyah and Singer; for example, it is injective on the copies of
$\mathbb{Z}$ in $\mathbb{B}_G$. There is one essential difference, however:
the $G$-signature given by $\theta$ can also detect real
representations (see also remark \ref{S_3}).
\begin{prop}
The map $\chi$ is an isomorphism whenever
$G\cong\mathbb{Z}/n\mathbb{Z}$, $n$ a positive integer.
\end{prop}
\begin{proof}
Conner and Floyd prove in their book \cite{CoFl} (see also
\cite[proposition 3.1]{EdEw}) that the cobordism
group of free $\mathbb{Z}/n\mathbb{Z}$-actions is isomorphic to
$H_2(\mathbb{Z}/n\mathbb{Z};\mathbb{Z})$ which is zero for every
$n$. Thus every free $G$-action is zero $G$-equivariant cobordant and
by theorem \ref{chi} $\chi$ is an isomorphism.
\end{proof}
As was already mentioned earlier, in \cite{Gr3} the author 
gives a complete computation of the group $\mathbb{B}_G$. This means we
can compute the $G$-equivariant cobordism group in dimension two, up
to the classes of free actions, for any finite group $G$.

\end{document}